\documentclass{article}
\usepackage{amsmath,amssymb,latexsym,wasysym,marvosym,textcomp,stmaryrd,tikz,multirow,mathtools,tkz-euclide,url,enumitem,amsthm,titlesec,caption,hyperref}
\usepackage[none]{hyphenat}
\usetikzlibrary{arrows,angles}

\title{On the NP-Completeness of Satisfying Certain Path and Loop Puzzles}
\author{Hadyn Tang}
\newcommand{\makeeye}[4]{
	\draw [color=#3] (#1+.9,#2+.45) arc (45:135:.32^.5);
	\draw [color=#3] (#1+.1,#2+.45) arc (225:315:.32^.5);
	\draw [color=#3] (#1+.5,#2+.45) circle (.32^.5-.4);
	\draw (#1+.25,#2+.75) node {\tiny #4};
}
\newcommand{\makeeyet}[4]{
	\draw [color=#3] (#1+.9,#2+.45) arc (45:135:.32^.5);
	\draw [color=#3] (#1+.1,#2+.45) arc (225:315:.32^.5);
	\draw [color=#3] (#1+.5,#2+.45) circle (.32^.5-.4);
	\draw (#1+.25,#2+.75) node {\scalebox{.3}{#4}};
}
\newcommand{\vtwist}[2]{
	\draw [color=white,line width=2pt] (#1-1.5,#2-1)--(#1+1.5,#2-1);
	\draw [color=white,line width=2pt] (#1-1.5,#2+1)--(#1+1.5,#2+1);
	\draw (#1,#2-1.5)--(#1,#2+1.5);
	\draw [color=blue,line width=1pt] (#1-1.5,#2)--(#1+1.5,#2);
	\draw [color=brown,line width=1pt] (#1-1.5,#2-1)--(#1-1,#2-1)--(#1-1,#2+1)--(#1-1.5,#2+1);
	\draw [color=brown,line width=1pt] (#1+1.5,#2-1)--(#1+1,#2-1)--(#1+1,#2+1)--(#1+1.5,#2+1);
}
\newcommand{\htwist}[2]{
	\draw [color=white,line width=2pt] (#1-1,#2-1.5)--(#1-1,#2+1.5);
	\draw [color=white,line width=2pt] (#1+1,#2-1.5)--(#1+1,#2+1.5);
	\draw (#1-1.5,#2)--(#1+1.5,#2);
	\draw [color=blue,line width=1pt] (#1,#2-1.5)--(#1,#2+1.5);
	\draw [color=brown,line width=1pt] (#1-1,#2-1.5)--(#1-1,#2-1)--(#1+1,#2-1)--(#1+1,#2-1.5);
	\draw [color=brown,line width=1pt] (#1-1,#2+1.5)--(#1-1,#2+1)--(#1+1,#2+1)--(#1+1,#2+1.5);
}
\begin{document}
	\maketitle
	
	\begin{center}
	\end{center}

	\begin{abstract}		
		\noindent ``Eye-Witless", ``Haisu" and ``Oriental House" are genres of logic puzzles invented by William Hu, and ``Detour" is a genre of logic puzzle invented by online user Guowen Zhang. Each of these puzzles revolves around constructing a path or loop through the cells of a grid according to certain constraints given by clues in the grid. We prove that deciding whether a particular puzzle in each of these genres is solvable is NP-complete.
	\end{abstract}

	\section{Introduction}
	
	We will first give an abridged explanation of the rules of each of the four puzzle genres, as well as an example puzzle and its solution. The word \textit{orthogonal} is here used to denote ``in the four cardinal directions". 
	
	\subsection {Eye-Witless}
	
	In Eye-Witless \cite{puzll} (fig 1), the player is given a rectangular grid, and the aim is to draw a loop passing orthogonally through the centres of its cells. Some eyes have been given with different colours: all eyes must be passed through by the loop. Additionally, if we consider the loop as a corridor with width one cell, same-coloured eyes can see the same distances, at the same relative directions. For example, if at a blue eye the loop turns, extending 3 cells down and 2 cells left before turning, at every blue eye the loop must go 3 cells before turning on the clockwise side of the bend, and 2 cells on the other side, though it might be 3 cells right and 2 cells up, 3 cells up and 2 cells left, or 3 cells left and 2 cells down instead. 
	
	\begin{figure}
		\centering
		\begin{tikzpicture}[font=\sffamily,scale=0.5]
		\foreach \x in {0,...,3}
		{
			\draw (0+6*\x,1)--(5+6*\x,1);
			\draw (0+6*\x,2)--(5+6*\x,2);
			\draw (0+6*\x,3)--(5+6*\x,3);
			\draw (0+6*\x,4)--(5+6*\x,4);
			\draw (1+6*\x,0)--(1+6*\x,5);
			\draw (2+6*\x,0)--(2+6*\x,5);
			\draw (3+6*\x,0)--(3+6*\x,5);
			\draw (4+6*\x,0)--(4+6*\x,5);
			\draw [line width=2pt] (0+6*\x,0)--(0+6*\x,5)--(5+6*\x,5)--(5+6*\x,0)--(0+6*\x,0)--(0+6*\x,5);
			\makeeye{6*\x+0}{1}{blue}{A};
			\makeeye{6*\x+0}{2}{blue}{A};
			\makeeye{6*\x+1}{0}{blue}{A};
			\makeeye{6*\x+1}{1}{blue}{A};
			\makeeye{6*\x+2}{0}{blue}{A};
			\makeeye{6*\x+2}{1}{green}{B};
			\makeeye{6*\x+2}{2}{green}{B};
			\makeeye{6*\x+2}{3}{green}{B};
			\makeeye{6*\x+3}{1}{green}{B};
			\makeeye{6*\x+4}{1}{brown}{C};
			\makeeye{6*\x+4}{2}{brown}{C};
		}
		\draw [color=red,line width=2pt] (7.5,2.5)--(6.5,2.5)--(6.5,1.5)--(7.5,1.5)--(7.5,0.5)--(8.5,0.5)--(8.5,1.5);
		\draw [color=red,line width=2pt] (13.5,3.5)--(13.5,2.5)--(12.5,2.5)--(12.5,1.5)--(13.5,1.5)--(13.5,0.5)--(14.5,0.5)--(14.5,1.5)--(15.5,1.5);
		\draw [color=red,line width=2pt] (14.5,2.5)--(14.5,3.5);
		\draw [color=red,line width=2pt] (16.5,0.5)--(16.5,2.5);
		\draw [color=red,line width=2pt] (18.5,1.5)--(19.5,1.5)--(19.5,0.5)--(20.5,0.5)--(20.5,1.5)--(21.5,1.5)--(21.5,0.5)--(22.5,0.5)--(22.5,3.5)--(21.5,3.5)--(21.5,2.5)--(20.5,2.5)--(20.5,3.5)
		--(19.5,3.5)--(19.5,2.5)--(18.5,2.5)--(18.5,1.5)--(19.5,1.5);
		\end{tikzpicture}
		\captionof{figure}{Eye-Witless example puzzle and steps to solution}
	\end{figure}
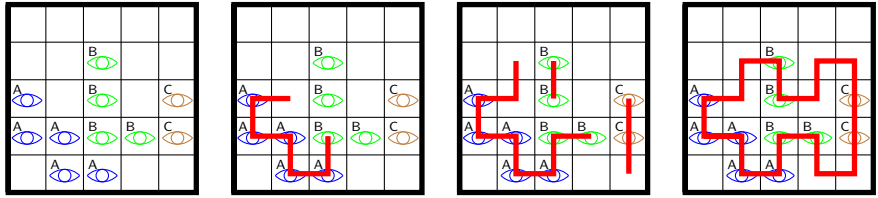
	
	\subsection{Haisu}
	
	In Haisu \cite{puzhaisu} (fig 2), the player is given a rectangular grid divided into regions, including two one-cell regions labelled with S for Start and F for Finish. The aim is to draw a directed orthogonal path from S to F passing through every cell exactly once. Some numbers have been given in the cells. A path may enter a region once or multiple times, if this is the $n^\text{th}$ time the path enters the region, it may only pass over black cells or instances of the number $n$.
	
	\begin{figure}
		\centering
		\begin{tikzpicture}[font=\sffamily,scale=0.5]
			\foreach \x in {0,...,3}
			{
			\draw (0+6*\x,1)--(5+6*\x,1);
			\draw (0+6*\x,2)--(5+6*\x,2);
			\draw (0+6*\x,3)--(5+6*\x,3);
			\draw (0+6*\x,4)--(5+6*\x,4);
			\draw (1+6*\x,0)--(1+6*\x,5);
			\draw (2+6*\x,0)--(2+6*\x,5);
			\draw (3+6*\x,0)--(3+6*\x,5);
			\draw (4+6*\x,0)--(4+6*\x,5);
			\draw [line width=2pt] (0+6*\x,0)--(0+6*\x,5)--(5+6*\x,5)--(5+6*\x,0)--(0+6*\x,0)--(0+6*\x,5);
			\draw [line width=2pt] (0+6*\x,3)--(1+6*\x,3)--(1+6*\x,2)--(2+6*\x,2)--(2+6*\x,1)--(3+6*\x,1)--(3+6*\x,0);
			\draw [line width=2pt] (2+6*\x,2)--(2+6*\x,4)--(3+6*\x,4)--(3+6*\x,2)--(4+6*\x,2)--(4+6*\x,1)--(3+6*\x,1);
			\draw [line width=2pt] (0+6*\x,1)--(1+6*\x,1)--(1+6*\x,0);
			\draw [line width=2pt] (5+6*\x,4)--(4+6*\x,4)--(4+6*\x,5);
			\draw (0.5+6*\x,0.5) node {\textbf S};
			\draw (4.5+6*\x,4.5) node {\textbf F};
			\draw (1.5+6*\x,1.5) node {2};
			\draw (1.5+6*\x,3.5) node {2};
			\draw (3.5+6*\x,1.5) node {4};
			\draw (3.5+6*\x,3.5) node {4};
			}
			\draw [color=red,line width=2pt] (6.5,3.5)--(6.5,4.5)--(7.5,4.5);
			\draw [color=red,line width=2pt] (9.5,0.5)--(10.5,0.5)--(10.5,1.5);
			\draw [color=red,line width=2pt] (7.5,2.5)--(9.5,2.5);
			\draw [color=red,line width=2pt] (7.5,1.5)--(8.5,1.5)--(8.5,0.5);
			\draw [color=red,line width=2pt] (12.5,0.5)--(14.5,0.5)--(14.5,1.5)--(13.5,1.5);
			\draw [color=red,line width=2pt] (12.5,3.5)--(12.5,4.5)--(13.5,4.5);
			\draw [color=red,line width=2pt] (13.5,2.5)--(15.5,2.5)--(15.5,0.5)--(16.5,0.5)--(16.5,3.5);
			\draw [color=red,line width=2pt] (18.5,0.5)--(20.5,0.5)--(20.5,1.5)--(18.5,1.5)--(18.5,4.5)--(20.5,4.5)--(20.5,3.5)--(19.5,3.5)--(19.5,2.5)--(21.5,2.5)--(21.5,0.5)--(22.5,0.5)--(22.5,3.5)--(21.5,3.5)--(21.5,4.5)--(22.5,4.5);
		\end{tikzpicture}
		\captionof{figure}{Haisu example puzzle and steps to solution}
	\end{figure}
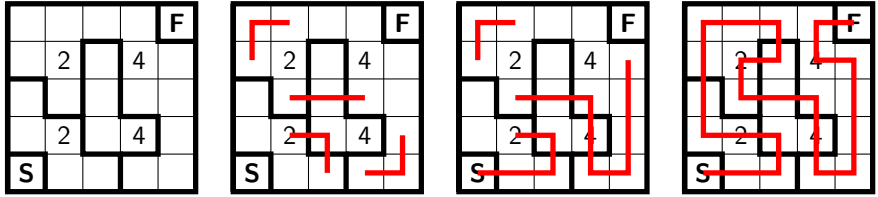
	
	\subsection{Oriental House}
	
	In Oriental House \cite{puzoh} (fig 3), similarly to Haisu, the player is given a rectangular grid divided into regions, including two one-cell regions labelled with S for Start and F for Finish. Again, the aim is to draw a directed orthogonal path from S to F passing through every cell exactly once. Some orthogonal arrows have been given in some cells: if the path passes over an arrow on a particular visit to a region, either it entered the region in that direction or will exit the region in that direction on this visit.
	
	\begin{figure}
		\centering
		\begin{tikzpicture}[font=\sffamily,scale=0.5]
			\foreach \x in {0,...,3}
			{
				\draw (0+6*\x,1)--(5+6*\x,1);
				\draw (0+6*\x,2)--(5+6*\x,2);
				\draw (0+6*\x,3)--(5+6*\x,3);
				\draw (0+6*\x,4)--(5+6*\x,4);
				\draw (1+6*\x,0)--(1+6*\x,5);
				\draw (2+6*\x,0)--(2+6*\x,5);
				\draw (3+6*\x,0)--(3+6*\x,5);
				\draw (4+6*\x,0)--(4+6*\x,5);
				\draw [line width=2pt] (0+6*\x,0)--(0+6*\x,5)--(5+6*\x,5)--(5+6*\x,0)--(0+6*\x,0)--(0+6*\x,5);
				\draw [line width=2pt] (0+6*\x,3)--(1+6*\x,3)--(1+6*\x,2)--(2+6*\x,2)--(2+6*\x,1)--(3+6*\x,1)--(3+6*\x,0);
				\draw [line width=2pt] (2+6*\x,2)--(2+6*\x,4)--(3+6*\x,4)--(3+6*\x,2)--(4+6*\x,2)--(4+6*\x,1)--(3+6*\x,1);
				\draw [line width=2pt] (0+6*\x,1)--(1+6*\x,1)--(1+6*\x,0);
				\draw [line width=2pt] (5+6*\x,4)--(4+6*\x,4)--(4+6*\x,5);
				\draw (0.5+6*\x,0.5) node {\textbf S};
				\draw (4.5+6*\x,4.5) node {\textbf F};
				\draw (1.5+6*\x,1.42) node {$\leftarrow$};
				\draw (1.5+6*\x,3.5) node {$\downarrow$};
				\draw (3.5+6*\x,1.5) node {$\uparrow$};
				\draw (3.5+6*\x,3.43) node {$\rightarrow$};
			}
			\draw [color=red,line width=2pt] (6.5,3.5)--(6.5,4.5)--(7.5,4.5);
			\draw [color=red,line width=2pt] (7.5,0.5)--(8.5,0.5);
			\draw [color=red,line width=2pt] (7.5,1.5)--(8.5,1.5);
			\draw [color=red,line width=2pt] (9.5,0.5)--(10.5,0.5)--(10.5,1.5);
			\draw [color=red,line width=2pt] (13.5,3.5)--(12.5,3.5)--(12.5,4.5)--(14.5,4.5)--(14.5,3.5);
			\draw [color=red,line width=2pt] (14.5,2.5)--(14.5,1.5)--(13.5,1.5);
			\draw [color=red,line width=2pt] (13.5,0.5)--(16.5,0.5)--(16.5,1.5)--(15.5,1.5)--(15.5,2.5)--(16.5,2.5)--(16.5,3.5);
			\draw [color=red,line width=2pt] (18.5,0.5)--(18.5,2.5)--(19.5,2.5)--(19.5,3.5)--(18.5,3.5)--(18.5,4.5)--(20.5,4.5)--(20.5,1.5)--(19.5,1.5)--(19.5,0.5)--(22.5,0.5)--(22.5,1.5)--(21.5,1.5)--(21.5,2.5)--(22.5,2.5)--(22.5,3.5)--(21.5,3.5)--(21.5,4.5)--(22.5,4.5);
		\end{tikzpicture}
		\captionof{figure}{Oriental House example puzzle and steps to solution}
	\end{figure}
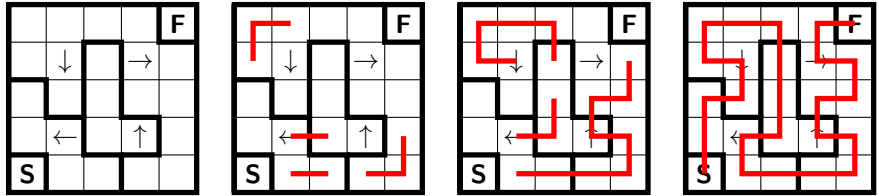
	
	\subsection{Detour}
	In Detour \cite{lmidetour} (fig 4), the player is given a rectangular grid divided into regions. The aim is to draw an orthogonal loop passing through each cell exactly once. Some regions may contain a small number: this number indicates the number of times the loop turns within this region.
	
	\begin{figure}
		\centering
		\begin{tikzpicture}[font=\sffamily,scale=0.4]
			\foreach \x in {0,...,3}
			{
				\draw (0+7*\x,1)--(6+7*\x,1);
				\draw (0+7*\x,2)--(6+7*\x,2);
				\draw (0+7*\x,3)--(6+7*\x,3);
				\draw (0+7*\x,4)--(6+7*\x,4);
				\draw (1+7*\x,0)--(1+7*\x,5);
				\draw (2+7*\x,0)--(2+7*\x,5);
				\draw (3+7*\x,0)--(3+7*\x,5);
				\draw (4+7*\x,0)--(4+7*\x,5);
				\draw (5+7*\x,0)--(5+7*\x,5);
				\draw [line width=2pt] (0+7*\x,0)--(0+7*\x,5)--(6+7*\x,5)--(6+7*\x,0)--(0+7*\x,0)--(0+7*\x,5);
				\draw [line width=2pt] (0+7*\x,3)--(1+7*\x,3)--(1+7*\x,2)--(2+7*\x,2)--(2+7*\x,1)--(3+7*\x,1)--(3+7*\x,0);
				\draw [line width=2pt] (2+7*\x,2)--(2+7*\x,4)--(3+7*\x,4)--(3+7*\x,2)--(4+7*\x,2)--(4+7*\x,1)--(3+7*\x,1);
				\draw [line width=2pt] (0+7*\x,1)--(1+7*\x,1)--(1+7*\x,0);
				\draw [line width=2pt] (5+7*\x,1)--(4+7*\x,1)--(4+7*\x,3)--(5+7*\x,3)--(5+7*\x,1)--(4+7*\x,1);
				\draw [line width=2pt] (4+7*\x,5)--(4+7*\x,4)--(5+7*\x,4)--(5+7*\x,5);
				\draw (0.3+7*\x,2.6) node {\tiny 1};
				\draw (2.3+7*\x,3.6) node {\tiny 1};
				\draw (4.3+7*\x,2.6) node {\tiny 0};
				\draw (4.3+7*\x,4.6) node {\tiny 1};
			}
			\draw [color=red,line width=2pt] (10.5,0.5)--(7.5,0.5)--(7.5,4.5)--(8.5,4.5);
			\draw [color=red,line width=2pt] (9.5,1.5)--(8.5,1.5)--(8.5,2.5);
			\draw [color=red,line width=2pt] (11.5,0.5)--(12.5,0.5)--(12.5,1.5);
			\draw [color=red,line width=2pt] (12.5,2.5)--(12.5,4.5)--(11.5,4.5)--(11.5,3.5);
			\draw [color=red,line width=2pt] (15.5,4.5)--(14.5,4.5)--(14.5,0.5)--(19.5,0.5)--(19.5,1.5)--(17.5,1.5);
			\draw [color=red,line width=2pt] (16.5,1.5)--(15.5,1.5)--(15.5,2.5);
			\draw [color=red,line width=2pt] (17.5,2.5)--(19.5,2.5)--(19.5,4.5)--(18.5,4.5)--(18.5,3.5);
			\draw [color=red,line width=2pt] (17.5,3.5)--(17.5,4.5)--(16.5,4.5);
			\draw [color=red,line width=2pt] (21.5,0.5)--(26.5,0.5)--(26.5,1.5)--(22.5,1.5)--(22.5,2.5)--(26.5,2.5)--(26.5,4.5)--(25.5,4.5)--(25.5,3.5)--(24.5,3.5)--(24.5,4.5)--(23.5,4.5)--(23.5,3.5)--(22.5,3.5)--(22.5,4.5)--(21.5,4.5)--(21.5,0.5)--(26.5,0.5);
		\end{tikzpicture}
		\captionof{figure}{Detour example puzzle and steps to solution}
	\end{figure}
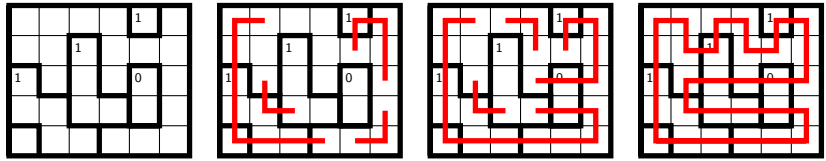

	\section{Membership in NP}
	
	Given a path or loop to the appropriate puzzle, it is easy to check whether it satisfies the conditions in P-time.
	
	Suppose the grid has dimensions $n\times m$, with $n\ge m$. Firstly, we can check that the path/loop passes over all the squares (or all the eyes, for Eye-Witless) in $O(n^2)$ time. Then for Haisu and Oriental House, we can pass over the path and check the number/arrow clues in $O(n^2)$ time. In Eye-Witless, we can find the length of the segments from a particular eye and their relative angles in $O(n)$ time, and iterating over the eyes takes $O(n^2)$, thus we can check all eyes are satisfied in $O(n^3)$ time. For Detour, we can iterate over the regions, counting the number of turns in each, in $O(n^2)$ time.
	
	\section{Reductions}
	
	\subsection{Eye-Witless: Cubic Planar Hamiltonicity}
	
	Motivated by the result established in \cite{pearlnp}, that we can reduce from cubic planar hamiltonicity to Masyu satisfiability, we aim to establish the same for the satisfiability of Eye-Witless in a similar manner. As observed on that page, the problem of cubic planar hamiltonicity is NP-complete (for example, it is a superset of the results in \cite{bicubic}). The first observation required is that we can draw a rectilinear realisation of any cubic planar graph (fig 5). 
	
	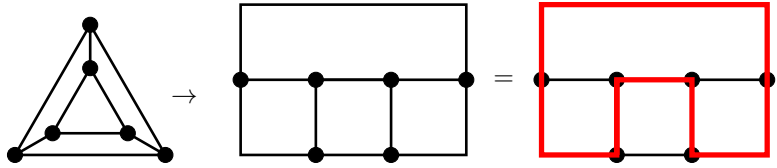
\begin{figure}
		\centering
		\begin{tikzpicture}[scale=.5]
		\draw [line width=1pt] (2,2*3^.5)--(0,0)--(4,0)--(2,2*3^.5)--(2,4/3^.5)--(1,1/3^.5)--(3,1/3^.5)--(2,4/3^.5)
		;
		\draw [line width=1pt] (0,0)--(1,1/3^.5);
		\draw [line width=1pt] (3,1/3^.5)--(4,0);
		\draw [fill=black] (0,0) circle (.2);
		\draw [fill=black] (4,0) circle (.2);
		\draw [fill=black] (2,2*3^.5) circle (.2);
		\draw [fill=black] (2,4/3^.5) circle (.2);
		\draw [fill=black] (1,1/3^.5) circle (.2);
		\draw [fill=black] (3,1/3^.5) circle (.2);
		\draw (4.5,1.5) node {$\rightarrow$};
		\draw [line width=1pt] (8,0)--(6,0)--(6,4)--(12,4)--(12,0)--(8,0)--(8,2)--(12,2);
		\draw [line width=1pt] (6,2)--(10,2)--(10,0);
		\draw [fill=black] (8,0) circle (.2);
		\draw [fill=black] (10,0) circle (.2);
		\draw [fill=black] (6,2) circle (.2);
		\draw [fill=black] (8,2) circle (.2);
		\draw [fill=black] (10,2) circle (.2);
		\draw [fill=black] (12,2) circle (.2);
		\draw (13,2) node {=};
		\draw [line width=1pt] (16,0)--(14,0)--(14,4)--(20,4)--(20,0)--(16,0)--(16,2)--(20,2);
		\draw [line width=1pt] (14,2)--(18,2)--(18,0);
		\draw [fill=black] (16,0) circle (.2);
		\draw [fill=black] (18,0) circle (.2);
		\draw [fill=black] (14,2) circle (.2);
		\draw [fill=black] (16,2) circle (.2);
		\draw [fill=black] (18,2) circle (.2);
		\draw [fill=black] (20,2) circle (.2);
		\draw [color=red,line width=2pt] (16,2)--(16,0)--(14,0)--(14,4)--(20,4)--(20,0)--(18,0)--(18,2)--(16,2)--(16,0);
		\end{tikzpicture}
		\captionof{figure}{Rectangular realisation of a graph, taken from \cite{pearlnp}}
	\end{figure}
	
	\cite{pearlnp} constructs the desired result by blocking off cells with a ``wall unit" that is constructed from several white pearls. Though we can't construct the equivalent white pearl in Eye-Witless, we can make something similar by use of a corner, allowing us to construct an eye through which the loop must go straight but turn on either side (fig 6). When chained, this eye type must form a wall.

	We can chain these eyes to form walls as in fig 6 (due to the turn on both sides nature of the A-eyes, we only need make these one unit thick), and then using these make wall units with one-cell-wide gaps between them as described in \cite{pearlnp}, with openings at the graph's vertices. Since the graph formed by the passageways between wall units is cubic, the loop can only visit each vertex at most once when we ignore the extra diversions caused by the wall units (since it can't pass an edge twice). Then we can place an extra A-eye at each vertex of the graph so that the loop passes each vertex once.

	\begin{figure}
		\centering
		\begin{tikzpicture}[font=\sffamily,scale=0.5]
			\draw (-1,1)--(2,1);
			\draw (-1,2)--(2,2);
			\draw (0,0)--(0,3);
			\draw (1,0)--(1,3);
			\draw [line width=2pt] (-1,0)--(-1,3)--(2,3);
			\makeeye{-1}{2}{green}{B};
			\makeeye{1}{2}{green}{B};
			\makeeye{0}{2}{blue}{A};
			\makeeye{-1}{1}{blue}{A};
			\draw (3,1.43) node {$\rightarrow$};
			\draw (4,1)--(7,1);
			\draw (4,2)--(7,2);
			\draw (5,0)--(5,3);
			\draw (6,0)--(6,3);
			\draw [line width=2pt] (4,0)--(4,3)--(7,3);
			\makeeye{4}{2}{green}{B};
			\makeeye{6}{2}{green}{B};
			\makeeye{5}{2}{blue}{A};
			\makeeye{4}{1}{blue}{A};
			\draw [color=red,line width=2pt] (5.5,0.5)--(4.5,0.5)--(4.5,2.5)--(6.5,2.5)--(6.5,0.5);
			\draw (8,1.43) node {$\rightarrow$};
			\draw (9,1)--(14,1);
			\draw (9,2)--(14,2);
			\draw (10,0)--(10,3);
			\draw (11,0)--(11,3);
			\draw (12,0)--(12,3);
			\draw (13,0)--(13,3);
			\makeeye{10}{1}{blue}{A};
			\makeeye{11}{1}{blue}{A};
			\makeeye{12}{1}{blue}{A};
			\draw [color=red,line width=2pt] (9.5,0.5)--(10.5,0.5)--(10.5,2.5)--(11.5,2.5)--(11.5,0.5)--(12.5,0.5)--(12.5,2.5)--(13.5,2.5);
			\draw (15,1.5) node {OR};
			\draw (16,1)--(21,1);
			\draw (16,2)--(21,2);
			\draw (17,0)--(17,3);
			\draw (18,0)--(18,3);
			\draw (19,0)--(19,3);
			\draw (20,0)--(20,3);
			\makeeye{17}{1}{blue}{A};
			\makeeye{18}{1}{blue}{A};
			\makeeye{19}{1}{blue}{A};
			\draw [color=red,line width=2pt] (16.5,2.5)--(17.5,2.5)--(17.5,0.5)--(18.5,0.5)--(18.5,2.5)--(19.5,2.5)--(19.5,0.5)--(20.5,0.5);
		\end{tikzpicture}
		\captionof{figure}{Forcing an eye type (A), and making a wall}
	\end{figure}
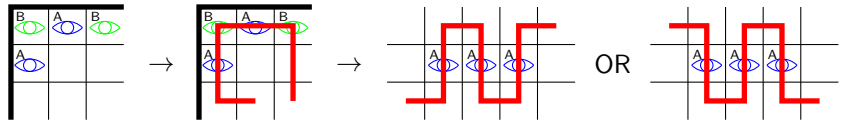

	There are two slight deviations. One is that the the corner configuration must appear: this can be achieved by having the outside wall open up to a corner at which to place the configuration, and the other is that sometimes due to parity we may need to extend the wall's length by a cell (which can be done by using B-eyes). Both are shown in figure 7, which provides a construction for the graph in figure 5. We note that, as stated in \cite{pearlnp}, the graph's representation can be done in $O(n^2)$ area, so the mapping from graphs to puzzles is polynomial and Eye-Witless satisfiability is NP-complete.

	\begin{figure}
		\centering
		\begin{tikzpicture}[font=\sffamily,scale=0.3]
			\foreach \x in {1,...,33} \draw (\x,0)--(\x,25);
			\foreach \x in {1,...,24} \draw (0,\x)--(34,\x);
			\foreach \x in {3,...,22}
			{
				\makeeyet{\x}{23}{blue}{A};
				\makeeyet{\x}{1}{blue}{A};
			}
			\foreach \x in {25,...,31}
			{
				\makeeyet{\x}{23}{blue}{A};
				\makeeyet{\x}{1}{blue}{A};
			}
			\foreach \x in {2,...,8}
			{
				\makeeyet{1}{\x}{blue}{A};
				\makeeyet{32}{\x}{blue}{A};
			}
			\foreach \x in {10,...,20}
			{
				\makeeyet{1}{\x}{blue}{A};
				\makeeyet{32}{\x}{blue}{A};
			}
			\makeeyet{2}{1}{blue}{A};
			\makeeyet{1}{24}{blue}{A};
			\makeeyet{0}{23}{blue}{A};
			\makeeyet{23}{23}{blue}{A};
			\makeeyet{32}{21}{blue}{A};
			\makeeyet{32}{22}{blue}{A};
			\makeeyet{0}{8}{green}{B};
			\makeeyet{0}{10}{green}{B};
			\makeeyet{33}{8}{green}{B};
			\makeeyet{33}{10}{green}{B};
			\makeeyet{0}{24}{green}{B};
			\makeeyet{2}{24}{green}{B};
			\makeeyet{23}{24}{green}{B};
			\makeeyet{25}{24}{green}{B};
			\draw [line width=2pt] (0,0)--(0,25)--(34,25)--(34,0)--(0,0)--(0,25);
			\foreach \x in {6,...,25}
			{
				\makeeyet{\x}{19}{blue}{A};
				\makeeyet{\x}{14}{blue}{A};
			}
			\makeeyet{26}{19}{blue}{A};
			\makeeyet{27}{19}{blue}{A};
			\makeeyet{5}{15}{blue}{A};
			\makeeyet{5}{16}{blue}{A};
			\makeeyet{5}{17}{blue}{A};
			\makeeyet{5}{18}{blue}{A};
			\makeeyet{28}{17}{blue}{A};
			\makeeyet{28}{18}{blue}{A};
			\foreach \x in {5,14,28}
			{
				\foreach \y in {6,...,9}
				{
					\makeeyet{\x}{\y}{blue}{A};
				}
			}
			\foreach \x in {6,15,24}
			{
				\foreach \xx in {0,...,3}
				{
					\makeeyet{\x+\xx}{10}{blue}{A};
				}
			}
			\makeeyet{6}{5}{blue}{A};
			\makeeyet{7}{5}{blue}{A};
			\makeeyet{15}{5}{blue}{A};
			\makeeyet{16}{5}{blue}{A};
			\makeeyet{26}{5}{blue}{A};
			\makeeyet{27}{5}{blue}{A};
			\makeeyet{10}{8}{blue}{A};
			\makeeyet{10}{9}{blue}{A};
			\makeeyet{19}{8}{blue}{A};
			\makeeyet{19}{9}{blue}{A};
			\makeeyet{23}{8}{blue}{A};
			\makeeyet{23}{9}{blue}{A};
			\makeeyet{3}{12}{blue}{A};
			\makeeyet{12}{12}{blue}{A};
			\makeeyet{21}{12}{blue}{A};
			\makeeyet{30}{12}{blue}{A};
			\makeeyet{12}{3}{blue}{A};
			\makeeyet{21}{3}{blue}{A};
		\end{tikzpicture}
		\begin{tikzpicture}[font=\sffamily,scale=0.3]
			\foreach \x in {1,...,33} \draw (\x,0)--(\x,25);
			\foreach \x in {1,...,24} \draw (0,\x)--(34,\x);
			\foreach \x in {3,...,22}
			{
				\makeeyet{\x}{23}{blue}{A};
				\makeeyet{\x}{1}{blue}{A};
			}
			\foreach \x in {25,...,31}
			{
				\makeeyet{\x}{23}{blue}{A};
				\makeeyet{\x}{1}{blue}{A};
			}
			\foreach \x in {2,...,8}
			{
				\makeeyet{1}{\x}{blue}{A};
				\makeeyet{32}{\x}{blue}{A};
			}
			\foreach \x in {10,...,20}
			{
				\makeeyet{1}{\x}{blue}{A};
				\makeeyet{32}{\x}{blue}{A};
			}
			\makeeyet{2}{1}{blue}{A};
			\makeeyet{1}{24}{blue}{A};
			\makeeyet{0}{23}{blue}{A};
			\makeeyet{23}{23}{blue}{A};
			\makeeyet{32}{21}{blue}{A};
			\makeeyet{32}{22}{blue}{A};
			\makeeyet{0}{8}{green}{B};
			\makeeyet{0}{10}{green}{B};
			\makeeyet{33}{8}{green}{B};
			\makeeyet{33}{10}{green}{B};
			\makeeyet{0}{24}{green}{B};
			\makeeyet{2}{24}{green}{B};
			\makeeyet{23}{24}{green}{B};
			\makeeyet{25}{24}{green}{B};
			\draw [line width=2pt] (0,0)--(0,25)--(34,25)--(34,0)--(0,0)--(0,25);
			\foreach \x in {6,...,25}
			{
				\makeeyet{\x}{19}{blue}{A};
				\makeeyet{\x}{14}{blue}{A};
			}
			\makeeyet{26}{19}{blue}{A};
			\makeeyet{27}{19}{blue}{A};
			\makeeyet{5}{15}{blue}{A};
			\makeeyet{5}{16}{blue}{A};
			\makeeyet{5}{17}{blue}{A};
			\makeeyet{5}{18}{blue}{A};
			\makeeyet{28}{17}{blue}{A};
			\makeeyet{28}{18}{blue}{A};
			\foreach \x in {5,14,28}
			{
				\foreach \y in {6,...,9}
				{
					\makeeyet{\x}{\y}{blue}{A};
				}
			}
			\foreach \x in {6,15,24}
			{
				\foreach \xx in {0,...,3}
				{
					\makeeyet{\x+\xx}{10}{blue}{A};
				}
			}
			\makeeyet{6}{5}{blue}{A};
			\makeeyet{7}{5}{blue}{A};
			\makeeyet{15}{5}{blue}{A};
			\makeeyet{16}{5}{blue}{A};
			\makeeyet{26}{5}{blue}{A};
			\makeeyet{27}{5}{blue}{A};
			\makeeyet{10}{8}{blue}{A};
			\makeeyet{10}{9}{blue}{A};
			\makeeyet{19}{8}{blue}{A};
			\makeeyet{19}{9}{blue}{A};
			\makeeyet{23}{8}{blue}{A};
			\makeeyet{23}{9}{blue}{A};
			\makeeyet{3}{12}{blue}{A};
			\makeeyet{12}{12}{blue}{A};
			\makeeyet{21}{12}{blue}{A};
			\makeeyet{30}{12}{blue}{A};
			\makeeyet{12}{3}{blue}{A};
			\makeeyet{21}{3}{blue}{A};
			
			\foreach \x in {2.5,4.5,...,20.5}
			{
				\draw [color=brown,line width=1pt] (\x,1.5)--(\x,0.5)--(\x+1,0.5)--(\x+1,2.5)--(\x+2,2.5)--(\x+2,1.5);
				\draw [color=brown,line width=1pt] (\x,23.5)--(\x,22.5)--(\x+1,22.5)--(\x+1,24.5)--(\x+2,24.5)--(\x+2,23.5);
			}
			\foreach \y in {3.5,5.5,10.5,12.5,...,18.5}
			{
				\draw [color=brown,line width=1pt] (1.5,\y)--(2.5,\y)--(2.5,\y+1)--(0.5,\y+1)--(0.5,\y+2)--(1.5,\y+2);
				\draw [color=brown,line width=1pt] (32.5,\y)--(31.5,\y)--(31.5,\y+1)--(33.5,\y+1)--(33.5,\y+2)--(32.5,\y+2);
			}
			\draw [color=brown,line width=1pt] (1.5,20.5)--(2.5,20.5)--(2.5,21.5)--(1.5,21.5)--(1.5,22.5)--(0.5,22.5)--(0.5,24.5)--(2.5,24.5)--(2.5,23.5);
			\draw [color=brown,line width=1pt] (1.5,7.5)--(2.5,7.5)--(2.5,8.5)--(0.5,8.5)--(0.5,10.5)--(1.5,10.5);
			\draw [color=brown,line width=1pt] (32.5,7.5)--(31.5,7.5)--(31.5,8.5)--(33.5,8.5)--(33.5,10.5)--(32.5,10.5);
			\draw [color=brown,line width=1pt] (1.5,3.5)--(0.5,3.5)--(0.5,2.5)--(2.5,2.5)--(2.5,1.5);
			\draw [color=brown,line width=1pt] (22.5,23.5)--(22.5,22.5)--(23.5,22.5)--(23.5,24.5)--(25.5,24.5)--(25.5,22.5)--(26.5,22.5)--(26.5,24.5)--(27.5,24.5)--(27.5,22.5)--(28.5,22.5)--(28.5,24.5)--(29.5,24.5)--(29.5,22.5)--(30.5,22.5)--(30.5,24.5)--(31.5,24.5)--(31.5,22.5)--(33.5,22.5)--(33.5,21.5)--(31.5,21.5)--(31.5,20.5)--(32.5,20.5);
			\foreach \x in {6.5,8.5,...,22.5}
			{
				\draw [color=brown,line width=1pt] (\x,14.5)--(\x,13.5)--(\x+1,13.5)--(\x+1,15.5)--(\x+2,15.5)--(\x+2,14.5);
				\draw [color=brown,line width=1pt] (\x,19.5)--(\x,20.5)--(\x+1,20.5)--(\x+1,18.5)--(\x+2,18.5)--(\x+2,19.5);
			}
			\draw [color=brown,line width=1pt] (6.5,14.5)--(6.5,15.5)--(4.5,15.5)--(4.5,16.5)--(6.5,16.5)--(6.5,17.5)--(4.5,17.5)--(4.5,18.5)--(6.5,18.5)--(6.5,19.5);
			\foreach \x in {4.5,13.5,22.5}
			{
				\draw [color=brown,line width=1pt] (\x+1,8.5)--(\x,8.5)--(\x,9.5)--(\x+2,9.5)--(\x+2,11.5)--(\x+3,11.5)--(\x+3,11.5)--(\x+3,9.5)--(\x+4,9.5)--(\x+4,11.5)--(\x+5,11.5)--(\x+5,9.5)--(\x+7,9.5)--(\x+7,8.5)--(\x+5,8.5)--(\x+5,8);
			}	
			\draw [color=brown,line width=1pt] (24.5,19.5)--(24.5,20.5)--(25.5,20.5)--(25.5,18.5)--(26.5,18.5)--(26.5,20.5)--(27.5,20.5)--(27.5,18.5)--(29.5,18.5)--(29.5,17.5)--(27.5,17.5)--(27.5,17);
			\draw [color=red,line width=2pt] (27.5,17)--(27.5,16.5)--(30.5,16.5)--(30.5,21.5)--(3.5,21.5)--(3.5,14.5)--(4.5,14.5)--(4.5,13.5)--(3.5,13.5)--(3.5,11.5)--(4.5,11.5)--(4.5,10.5)--(3.5,10.5)--(3.5,3.5)--(10.5,3.5)--(10.5,4.5)--(11.5,4.5)--(11.5,3.5)--(13.5,3.5)--(13.5,5.5)--(8.5,5.5)--(8.5,6.5)--(8,6.5);
			\draw [color=brown,line width=1pt] (8,6.5)--(7.5,6.5)--(7.5,4.5)--(6.5,4.5)--(6.5,6.5)--(4.5,6.5)--(4.5,7.5)--(6.5,7.5)--(6.5,8.5)--(5.5,8.5); 
			\draw [color=red,line width=2pt] (9.5,8)--(9.5,7.5)--(12.5,7.5)--(12.5,10.5)--(11.5,10.5)--(11.5,12.5)--(13.5,12.5)--(13.5,11.5)--(14.5,11.5)--(14.5,12.5)--(19.5,12.5)--(19.5,11.5)--(20.5,11.5)--(20.5,12.5)--(22.5,12.5)--(22.5,10.5)--(21.5,10.5)--(21.5,7.5)--(18.5,7.5)--(18.5,8);
			\draw [color=brown,line width=1pt] (14.5,8.5)--(15.5,8.5)--(15.5,7.5)--(13.5,7.5)--(13.5,6.5)--(15.5,6.5)--(15.5,4.5)--(16.5,4.5)--(16.5,6.5)--(17,6.5);
			\draw [color=red,line width=2pt] (17,6.5)--(20.5,6.5)--(20.5,3.5)--(22.5,3.5)--(22.5,4.5)--(23.5,4.5)--(23.5,0.5)--(23,0.5);
			\draw [color=brown,line width=1pt] (23,0.5)--(22.5,0.5)--(22.5,1.5);
			\draw [color=brown,line width=1pt] (23.5,8.5)--(24.5,8.5)--(24.5,8);
			\draw [color=red,line width=2pt] (24.5,8)--(24.5,0.5)--(25,0.5);
			\draw [color=brown,line width=1pt] (25,0.5)--(25.5,0.5)--(25.5,2.5)--(26.5,2.5)--(26.5,0.5)--(27.5,0.5)--(27.5,2.5)--(28.5,2.5)--(28.5,0.5)--(29.5,0.5)--(29.5,2.5)--(30.5,2.5)--(30.5,0.5)--(31.5,0.5)--(31.5,2.5)--(33.5,2.5)--(33.5,3.5)--(32.5,3.5);
			\draw [color=brown,line width=1pt] (27.5,8)--(27.5,7.5)--(29.5,7.5)--(29.5,6.5)--(27.5,6.5)--(27.5,4.5)--(26.5,4.5)--(26.5,6.5)--(26,6.5);
			\draw [color=red,line width=2pt] (26,6.5)--(25.5,6.5)--(25.5,3.5)--(30.5,3.5)--(30.5,10.5)--(29.5,10.5)--(29.5,11.5)--(30.5,11.5)--(30.5,13.5)--(26.5,13.5)--(26.5,15.5)--(26,15.5);
			\draw [color=brown,line width=1pt] (26,15.5)--(25.5,15.5)--(25.5,13.5)--(24.5,13.5)--(24.5,14.5);
		\end{tikzpicture}
		\captionof{figure}{Eye-Witless example construction with solution from fig 5}
	\end{figure}
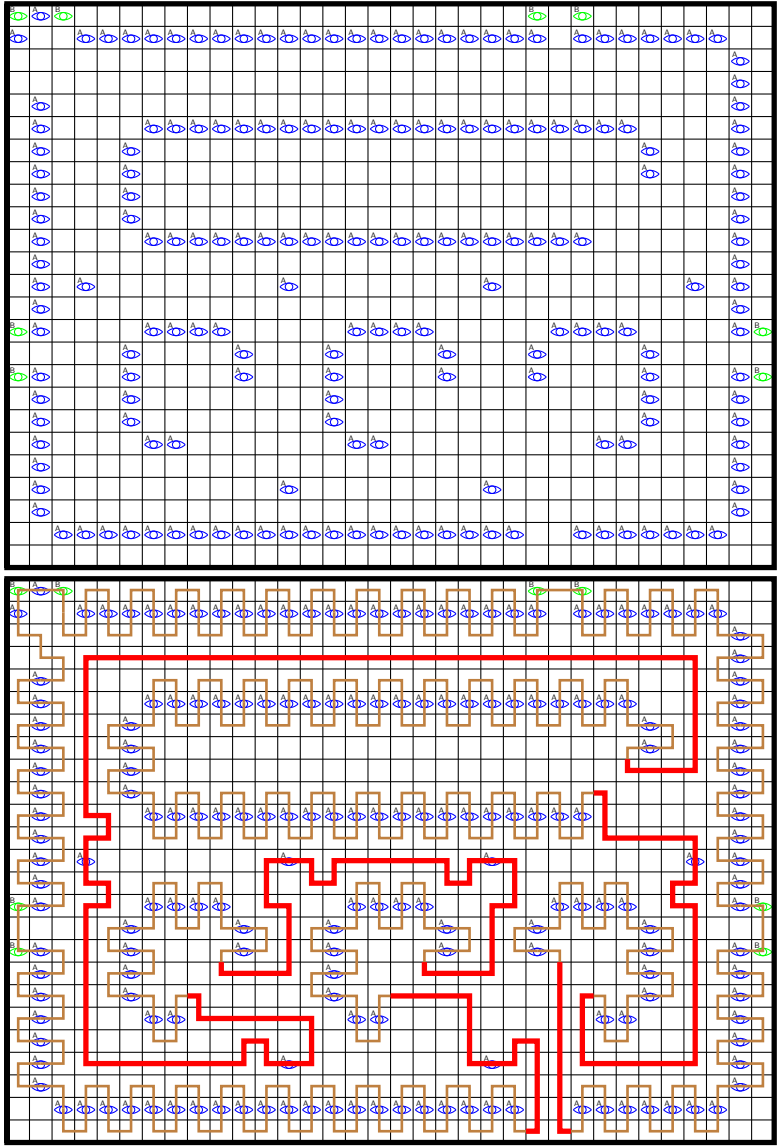
	
	\subsection{Haisu: Bipartite Cubic Planar Hamiltonicity}
	
	Unlike Eye-Witness, the other three genres require every cell be passed through. This brings with it two issues: needing to ensure that our design allows for every cell to be passed through, and the associated issues with parity. To help deal with parity, we resort to a stronger result found in \cite{bicubic}: that the Hamiltonicity of a given bipartite cubic planar graph is also NP-complete.
	
	Again, the first step is to draw the graph in a rectilinear fashion as we did for Eye-Witness. Since the graph is bipartite, we may also colour its vertices black and white. Then, we impose the additional constraint that, if the rectilinear drawing of the grid is superimposed over a checkerboard, the colours of the respective vertices should match the cell that they are on (fig 8). This can be achieved by taking the rectangular realisation, dividing each of its into $3\times 3$ grids and shifting any non-matching colours one cell down, as shown in figure 9. We also need to ensure that the grid has an even dimension: we can do this by adding an extra empty column if necessary.
	
	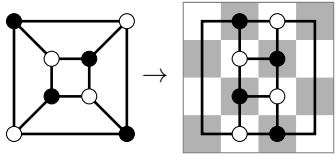
\begin{figure}
		\centering
		\begin{tikzpicture}[scale=0.5]
			\draw [line width=1pt] (0,0)--(0,3)--(3,3)--(3,0)--(0,0)--(1,1)--(1,2)--(2,2)--(2,1)--(1,1);
			\draw [line width=1pt] (1,2)--(0,3);
			\draw [line width=1pt] (2,1)--(3,0);
			\draw [line width=1pt] (2,2)--(3,3);
			\draw [fill=white] (3,3) circle (.2);
			\draw [fill=white] (1,2) circle (.2);
			\draw [fill=white] (2,1) circle (.2);
			\draw [fill=white] (0,0) circle (.2);
			\draw [fill=black] (0,3) circle (.2);
			\draw [fill=black] (2,2) circle (.2);
			\draw [fill=black] (1,1) circle (.2);
			\draw [fill=black] (3,0) circle (.2);
			\draw (3.75,1.5) node {$\rightarrow$};
			\fill [gray!60] (4.5,-0.5)--(5.5,-0.5)--(5.5,3.5)--(6.5,3.5)--(6.5,-0.5)--(7.5,-0.5)--(7.5,3.5)--(8.5,3.5)--(8.5,2.5)--(4.5,2.5)--(4.5,1.5)--(8.5,1.5)--(8.5,0.5)--(4.5,0.5)--cycle;
			\draw [color=gray] (4.5,3.5)--(8.5,3.5)--(8.5,-0.5)--(4.5,-0.5)--(4.5,3.5)--(8.5,3.5);
			\draw [line width=1pt] (6,0)--(6,3)--(5,3)--(5,0)--(8,0)--(8,3)--(7,3)--(7,0);
			\draw [line width=1pt] (6,1)--(7,1);
			\draw [line width=1pt] (6,2)--(7,2);
			\draw [line width=1pt] (6,3)--(7,3);
			\draw [fill=white] (7,3) circle (.2);
			\draw [fill=white] (6,2) circle (.2);
			\draw [fill=white] (7,1) circle (.2);
			\draw [fill=white] (6,0) circle (.2);
			\draw [fill=black] (6,3) circle (.2);
			\draw [fill=black] (7,2) circle (.2);
			\draw [fill=black] (6,1) circle (.2);
			\draw [fill=black] (7,0) circle (.2);
		\end{tikzpicture}
		\captionof{figure}{Appropriate bipartite cubic planar graph rectilinear realisation}
	\end{figure}

	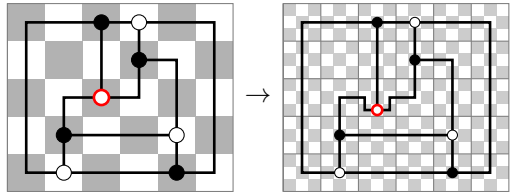
\begin{figure}
		\centering
		\begin{tikzpicture}[scale=1/6]
			\foreach \x in {-1.5,4.5,10.5}
			{
				\foreach \y in {-1.5,4.5,10.5} \fill [gray!60] (\x,\y)--(\x,\y+3)--(\x+3,\y+3)--(\x+3,\y)--cycle;
			}
			\foreach \x in {1.5,7.5,13.5}
			{
				\foreach \y in {1.5,7.5} \fill [gray!60] (\x,\y)--(\x,\y+3)--(\x+3,\y+3)--(\x+3,\y)--cycle;
			}
			\draw [color=gray] (-1.5,-1.5)--(-1.5,13.5)--(16.5,13.5)--(16.5,-1.5)--(-1.5,-1.5)--(-1.5,13.5);
			\draw [line width=1pt] (6,12)--(6,6)--(3,6)--(3,0)--(0,0)--(0,12)--(15,12)--(15,0)--(12,0)--(12,9)--(9,9)--(9,6)--(6,6);
			\draw [line width=1pt] (3,0)--(12,0);
			\draw [line width=1pt] (3,3)--(12,3);
			\draw [line width=1pt] (9,9)--(9,12);
			\draw [fill=white] (3,0) circle (.6);
			\draw [color=red,line width=1pt,fill=white] (6,6) circle (.6);
			\draw [fill=white] (9,12) circle (.6);
			\draw [fill=white] (12,3) circle (.6);
			\draw [fill=black] (3,3) circle (.6);
			\draw [fill=black] (6,12) circle (.6);
			\draw [fill=black] (9,9) circle (.6);
			\draw [fill=black] (12,0) circle (.6);
			\draw (18.5,6) node {$\rightarrow$};
			\foreach \x in {20.5,22.5,...,36.5}
			{
				\foreach \y in {-1.5,0.5,...,12.5} \fill [gray!40] (\x,\y)--(\x,\y+1)--(\x+1,\y+1)--(\x+1,\y)--cycle;
			}
			\foreach \x in {21.5,23.5,...,37.5}
			{
				\foreach \y in {-0.5,1.5,...,11.5} \fill [gray!40] (\x,\y)--(\x,\y+1)--(\x+1,\y+1)--(\x+1,\y)--cycle;
			}
			\foreach \x in {23.5,26.5,...,35.5} \draw [color=gray] (\x,-1.5)--(\x,13.5);
			\foreach \y in {1.5,4.5,7.5,10.5} \draw [color=gray] (20.5,\y)--(38.5,\y);
			\draw [color=gray] (20.5,-1.5)--(20.5,13.5)--(38.5,13.5)--(38.5,-1.5)--(20.5,-1.5)--(20.5,13.5);
			\draw [line width=1pt] (28,12)--(28,5)--(27,5)--(27,6)--(25,6)--(25,0)--(22,0)--(22,12)--(37,12)--(37,0)--(34,0)--(34,9)--(31,9)--(31,6)--(29,6)--(29,5)--(28,5);
			\draw [line width=1pt] (25,0)--(34,0);
			\draw [line width=1pt] (25,3)--(34,3);
			\draw [line width=1pt] (31,9)--(31,12);
			\draw [fill=white] (25,0) circle (.4);
			\draw [color=red,line width=1pt,fill=white] (28,5) circle (.4);
			\draw [fill=white] (31,12) circle (.4);
			\draw [fill=white] (34,3) circle (.4);
			\draw [fill=black] (25,3) circle (.4);
			\draw [fill=black] (28,12) circle (.4);
			\draw [fill=black] (31,9) circle (.4);
			\draw [fill=black] (34,0) circle (.4);
		\end{tikzpicture}
		\captionof{figure}{Fixing an incorrect realisation due to incorrect vertex}
	\end{figure}
	
	Once that has been completed, we next wish to change the graph so that all the blank cells are on one row. We do this by `unfolding' the columns onto a line in a snake-like manner, with connection between columns being raised alternately over and under the row, at an even height. Then, to fill empty cells not in the given row we exercise one of two options: either there is the outside of one connection to its left or right, so we can extrude a segment from the connection outwards (X); or it is closer to the row than some connection but further from the main row than any other contained connections, so we can extrude the horizontal section of the connection towards the main row (Y). This process is shown in figure 10 (with a simple multigraph because the process effectively squares the width and height).
	
	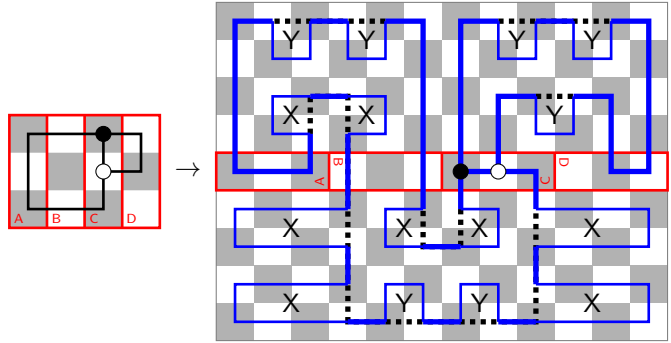
\begin{figure}
		\centering
		\begin{tikzpicture}[scale=0.25,font=\sffamily]
			\fill [gray!60] (0,-3)--(0,-1)--(8,-1)--(8,1)--(6,1)--(6,-3)--(4,-3)--(4,3)--(6,3)--(6,1)--(0,1)--(0,3)--(2,3)--(2,-3)--cycle;
			\draw [color=red,line width=1pt] (0,-3)--(8,-3)--(8,3)--(0,3)--cycle;
			\foreach \x in {2,4,6} \draw [color=red,line width=1pt] (\x,-3)--(\x,3);
			\draw [color=red] (0.5,-2.5) node {\tiny A};
			\draw [color=red] (2.5,-2.5) node {\tiny B};
			\draw [color=red] (4.5,-2.5) node {\tiny C};
			\draw [color=red] (6.5,-2.5) node {\tiny D};
			\draw [line width=1pt] (5,2)--(1,2)--(1,-2)--(5,-2)--(5,2)--(7,2)--(7,0)--(5,0);
			\draw [fill=white] (5,0) circle (.4);
			\draw [fill=black] (5,2) circle (.4);
			\draw (9.5,0) node {$\rightarrow$};
			\foreach \x in {11,15,...,31}
			{
				\foreach \y in {-9,-5,-1,3,7} \fill [gray!60] (\x,\y)--(\x+2,\y)--(\x+2,\y+2)--(\x,\y+2)--cycle;
			}
			\foreach \x in {13,17,...,33}
			{
				\foreach \y in {-7,-3,1,5} \fill [gray!60] (\x,\y)--(\x+2,\y)--(\x+2,\y+2)--(\x,\y+2)--cycle;
			}
			\draw [color=gray] (11,-9)--(11,9)--(35,9)--(35,-9)--cycle;
			\draw [color=red,line width=1pt] (11,-1)--(11,1)--(35,1)--(35,-1)--cycle;
			\foreach \x in {17,23,29} \draw [color=red,line width=1pt] (\x,-1)--(\x,1);
			\draw [dotted,line width=2pt] (24,-2)--(24,-4)--(22,-4)--(22,-2);
			\draw [dotted,line width=2pt] (20,8)--(14,8);
			\draw [dotted,line width=2pt] (16,2)--(16,4)--(18,4)--(18,-8)--(28,-8)--(28,-2);
			\draw [dotted,line width=2pt] (28,4)--(30,4);
			\draw [dotted,line width=2pt] (26,8)--(32,8);
			\draw [color=blue,line width=2pt] (24,0)--(24,-2);
			\draw [color=blue,line width=1pt] (24,-2)--(26,-2)--(26,-4)--(24,-4);
			\draw [color=blue,line width=2pt] (24,-4)--(22,-4);
			\draw [color=blue,line width=1pt] (22,-4)--(20,-4)--(20,-2)--(22,-2);
			\draw [color=blue,line width=2pt] (22,-2)--(22,8)--(20,8);
			\draw [color=blue,line width=1pt] (20,8)--(20,6)--(18,6)--(18,8);
			\draw [color=blue,line width=2pt] (18,8)--(16,8);
			\draw [color=blue,line width=1pt] (16,8)--(16,6)--(14,6)--(14,8);
			\draw [color=blue,line width=2pt] (14,8)--(12,8)--(12,0)--(16,0)--(16,2);
			\draw [color=blue,line width=1pt] (16,2)--(14,2)--(14,4)--(16,4);
			\draw [color=blue,line width=2pt] (16,4)--(18,4);
			\draw [color=blue,line width=1pt] (18,4)--(20,4)--(20,2)--(18,2);
			\draw [color=blue,line width=2pt] (18,2)--(18,-2);
			\draw [color=blue,line width=1pt] (18,-2)--(12,-2)--(12,-4)--(18,-4);
			\draw [color=blue,line width=2pt] (18,-4)--(18,-6);
			\draw [color=blue,line width=1pt] (18,-6)--(12,-6)--(12,-8)--(18,-8);
			\draw [color=blue,line width=2pt] (18,-8)--(20,-8);
			\draw [color=blue,line width=1pt] (20,-8)--(20,-6)--(22,-6)--(22,-8);
			\draw [color=blue,line width=2pt] (22,-8)--(24,-8);
			\draw [color=blue,line width=1pt] (24,-8)--(24,-6)--(26,-6)--(26,-8);
			\draw [color=blue,line width=2pt] (26,-8)--(28,-8);
			\draw [color=blue,line width=1pt] (28,-8)--(34,-8)--(34,-6)--(28,-6);
			\draw [color=blue,line width=2pt] (28,-6)--(28,-4);
			\draw [color=blue,line width=1pt] (28,-4)--(34,-4)--(34,-2)--(28,-2);
			\draw [color=blue,line width=2pt] (28,-2)--(28,0)--(24,0)--(24,8)--(26,8);
			\draw [color=blue,line width=1pt] (26,8)--(26,6)--(28,6)--(28,8);
			\draw [color=blue,line width=2pt] (28,8)--(30,8);
			\draw [color=blue,line width=1pt] (30,8)--(30,6)--(32,6)--(32,8);
			\draw [color=blue,line width=2pt] (32,8)--(34,8)--(34,0)--(32,0)--(32,4)--(30,4);
			\draw [color=blue,line width=1pt] (30,4)--(30,2)--(28,2)--(28,4);
			\draw [color=blue,line width=2pt] (28,4)--(26,4)--(26,0);
			\draw [fill=white] (26,0) circle (.4);
			\draw [fill=black] (24,0) circle (.4);
			\draw [color=red] (16.5,-0.5) node {\tiny \rotatebox{90}{A}};
			\draw [color=red] (17.5,0.5) node {\tiny \rotatebox{270}{B}};
			\draw [color=red] (28.5,-0.5) node {\tiny \rotatebox{90}{C}};
			\draw [color=red] (29.5,0.5) node {\tiny \rotatebox{270}{D}};
			\draw (15,-7) node {X};
			\draw (15,-3) node {X};
			\draw (15,3) node {X};
			\draw (15,7) node {Y};
			\draw (19,3) node {X};
			\draw (19,7) node {Y};
			\draw (21,-3) node {X};
			\draw (21,-7) node {Y};
			\draw (25,-3) node {X};
			\draw (25,-7) node {Y};
			\draw (27,7) node {Y};
			\draw (29,3) node {Y};
			\draw (31,7) node {Y};
			\draw (31,-3) node {X};
			\draw (31,-7) node {X};
		\end{tikzpicture}
		\captionof{figure}{Expanding the realisation into a row}
	\end{figure}

	For the next step, we split the expanded grid further into another subgrid a third of the size, with cells of the grid corresponding to monominoes (with the same colour as their corresponding cell), edges of the grid corresponding to dominoes (perpendicular to the edge) and corners of the grid corresponding to a $2\times 2$ horizontal domino squares (shown in the first diagram of figure 12).
	
	We would like all the monominoes to be part of the edges of the graph. However, some monominoes in the centre row correspond to empty squares. We notice that by considering the directed graph with edges going towards white vertices, the total outdegree is 3 times the number of black vertices, but this is also equal to the total indegree which is 3 times the number of white vertices, so there are the same number of each. So the vertices cover an equal number of white and black squares, and since all edges are also black-white they have odd length and also pass over an equal number of white and black squares. Thus, since the grid has an even dimension and thus an equal number of white and black squares, it has an equal number of empty white and black squares.
	
	Thus, as a remedy, we can pair unused black and white monominoes along the row so that there exist paths between them that don't cross, where the paths go from $2\times2$ cell to cell, and start and end at a cell with the monomino at its corner (for example by using a queue for excess monominoes of each colour). (We might need to expand the height of the grid, by adjusting the top connection up an appropriate even number of rows, but using the row-expanded construction, the height of the grid above the row required by this construction is at most the number of empty squares.) Then, we can switch the black monomino with an edge domino, repeatedly substitute the dominoes along the edges in the path and join the black and white monominoes at the end of the path, using the substitutions in figure 11 (and illustration on a different graph in figure 12).
	
	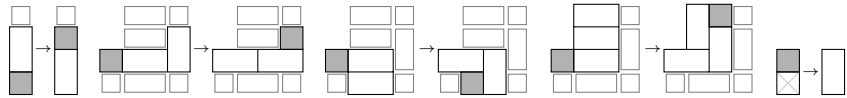
\begin{figure}
		\centering
		\begin{tikzpicture}[scale=0.3]
			\draw [fill=gray!60] (0,0)--(1,0)--(1,1)--(0,1)--cycle;
			\draw (0,1)--(0,3)--(1,3)--(1,1)--cycle;
			\draw (1.5,2) node {\tiny $\rightarrow$};
			\draw (2,0)--(2,2)--(3,2)--(3,0)--cycle;
			\draw [fill=gray!60] (2,2)--(2,3)--(3,3)--(3,2)--cycle;
			\draw [color=gray] (0.1,3.1)--(0.9,3.1)--(0.9,3.9)--(0.1,3.9)--cycle;
			\draw [color=gray] (2.1,3.1)--(2.9,3.1)--(2.9,3.9)--(2.1,3.9)--cycle;
			\draw [fill=gray!60] (4,1)--(4,2)--(5,2)--(5,1)--cycle;
			\draw (7,2)--(5,2)--(5,1)--(8,1)--(8,3)--(7,3)--(7,1);
			\draw [color=gray] (4.1,0.1)--(4.1,0.9)--(4.9,0.9)--(4.9,0.1)--cycle;
			\draw [color=gray] (5.1,0.1)--(5.1,0.9)--(6.9,0.9)--(6.9,0.1)--cycle;
			\draw [color=gray] (7.1,0.1)--(7.1,0.9)--(7.9,0.9)--(7.9,0.1)--cycle;
			\draw [color=gray] (5.1,3.1)--(5.1,3.9)--(6.9,3.9)--(6.9,3.1)--cycle;
			\draw [color=gray] (5.1,2.1)--(5.1,2.9)--(6.9,2.9)--(6.9,2.1)--cycle;
			\draw [color=gray] (7.1,3.1)--(7.1,3.9)--(7.9,3.9)--(7.9,3.1)--cycle;
			\draw [fill=gray!60] (12,2)--(12,3)--(13,3)--(13,2)--cycle;
			\draw (11,2)--(11,1)--(13,1)--(13,2)--(9,2)--(9,1)--(11,1);
			\draw [color=gray] (9.1,0.1)--(9.1,0.9)--(9.9,0.9)--(9.9,0.1)--cycle;
			\draw [color=gray] (10.1,0.1)--(10.1,0.9)--(11.9,0.9)--(11.9,0.1)--cycle;
			\draw [color=gray] (12.1,0.1)--(12.1,0.9)--(12.9,0.9)--(12.9,0.1)--cycle;
			\draw [color=gray] (10.1,3.1)--(10.1,3.9)--(11.9,3.9)--(11.9,3.1)--cycle;
			\draw [color=gray] (10.1,2.1)--(10.1,2.9)--(11.9,2.9)--(11.9,2.1)--cycle;
			\draw [color=gray] (12.1,3.1)--(12.1,3.9)--(12.9,3.9)--(12.9,3.1)--cycle;
			\draw (8.5,2) node {\tiny $\rightarrow$};
			\draw [fill=gray!60] (14,1)--(14,2)--(15,2)--(15,1)--cycle;
			\draw (15,1)--(17,1)--(17,0)--(15,0)--(15,2)--(17,2)--(17,1);
			\draw [color=gray] (14.1,0.1)--(14.1,0.9)--(14.9,0.9)--(14.9,0.1)--cycle;
			\draw [color=gray] (17.1,0.1)--(17.1,0.9)--(17.9,0.9)--(17.9,0.1)--cycle;
			\draw [color=gray] (17.1,3.1)--(17.1,3.9)--(17.9,3.9)--(17.9,3.1)--cycle;
			\draw [color=gray] (15.1,3.1)--(15.1,3.9)--(16.9,3.9)--(16.9,3.1)--cycle;
			\draw [color=gray] (15.1,2.1)--(15.1,2.9)--(16.9,2.9)--(16.9,2.1)--cycle;
			\draw [color=gray] (17.1,1.1)--(17.1,2.9)--(17.9,2.9)--(17.9,1.1)--cycle;
			\draw (18.5,2) node {\tiny $\rightarrow$};
			\draw [fill=gray!60] (20,0)--(21,0)--(21,1)--(20,1)--cycle;
			\draw (21,1)--(21,2)--(22,2)--(22,0)--(21,0)--(21,1)--(19,1)--(19,2)--(21,2);
			\draw [color=gray] (19.1,0.1)--(19.1,0.9)--(19.9,0.9)--(19.9,0.1)--cycle;
			\draw [color=gray] (22.1,0.1)--(22.1,0.9)--(22.9,0.9)--(22.9,0.1)--cycle;
			\draw [color=gray] (22.1,3.1)--(22.1,3.9)--(22.9,3.9)--(22.9,3.1)--cycle;
			\draw [color=gray] (20.1,3.1)--(20.1,3.9)--(21.9,3.9)--(21.9,3.1)--cycle;
			\draw [color=gray] (20.1,2.1)--(20.1,2.9)--(21.9,2.9)--(21.9,2.1)--cycle;
			\draw [color=gray] (22.1,1.1)--(22.1,2.9)--(22.9,2.9)--(22.9,1.1)--cycle;
			\draw [fill=gray!60] (24,1)--(24,2)--(25,2)--(25,1)--cycle;
			\draw (25,2)--(27,2)--(27,1)--(25,1)--(25,4)--(27,4)--(27,3)--(25,3);
			\draw (27,2)--(27,3);
			\draw [color=gray] (24.1,0.1)--(24.1,0.9)--(24.9,0.9)--(24.9,0.1)--cycle;
			\draw [color=gray] (27.1,0.1)--(27.1,0.9)--(27.9,0.9)--(27.9,0.1)--cycle;
			\draw [color=gray] (27.1,3.1)--(27.1,3.9)--(27.9,3.9)--(27.9,3.1)--cycle;
			\draw [color=gray] (25.1,0.1)--(25.1,0.9)--(26.9,0.9)--(26.9,0.1)--cycle;
			\draw [color=gray] (27.1,1.1)--(27.1,2.9)--(27.9,2.9)--(27.9,1.1)--cycle;
			\draw (28.5,2) node {\tiny $\rightarrow$};
			\draw [fill=gray!60] (31,3)--(32,3)--(32,4)--(31,4)--cycle;
			\draw (31,3)--(31,1)--(32,1)--(32,3)--(31,3)--(31,4)--(30,4)--(30,2)--(29,2)--(29,1)--(31,1);
			\draw (30,2)--(31,2);
			\draw [color=gray] (29.1,0.1)--(29.1,0.9)--(29.9,0.9)--(29.9,0.1)--cycle;
			\draw [color=gray] (32.1,0.1)--(32.1,0.9)--(32.9,0.9)--(32.9,0.1)--cycle;
			\draw [color=gray] (32.1,3.1)--(32.1,3.9)--(32.9,3.9)--(32.9,3.1)--cycle;
			\draw [color=gray] (30.1,0.1)--(30.1,0.9)--(31.9,0.9)--(31.9,0.1)--cycle;
			\draw [color=gray] (32.1,1.1)--(32.1,2.9)--(32.9,2.9)--(32.9,1.1)--cycle;
			\draw [fill=gray!60] (34,1)--(34,2)--(35,2)--(35,1)--cycle;
			\draw [color=gray!40] (34,0)--(35,1);
			\draw [color=gray!40] (34,1)--(35,0);
			\draw (34,1)--(34,0)--(35,0)--(35,1)--cycle;
			\draw (35.5,1) node {\tiny $\rightarrow$};
			\draw (36,0)--(36,2)--(37,2)--(37,0)--cycle;
		\end{tikzpicture}
		\captionof{figure}{Substitutions to remove monominoes}
	\end{figure}

	Now we have almost all the groundwork to make our construction. We need to divide the grid into one-monomino `white' regions containing the white cells and large `black' regions containing black cells, so that the adjacency graph between white and black regions corresponds to the original graph. We can do this by assigning edges in the graph between black and white vertices to the corresponding black vertex. This will guarantee all given adjacencies from the original graph are represented. And to ensure no extra adjacencies are included, we have to ensure that the fourth domino connected to any white monomino (if such a fourth domino exists) is included in one of the three regions connected to the vertex. We can easily do this because it's not automatically included in an edge by assumption, and thus we can just attach it, possibly with another domino, to an adjoining edge on either side of the monomino. (See fig 12, olive region.) We can then assign the other dominoes in any manner: this is shown in the last picture in figure 12: monominoes/dominoes in olive or with the blue graph passing over them are required to be part of the region they are in, but the other dominoes could have been in another region if there were another white one.
	
	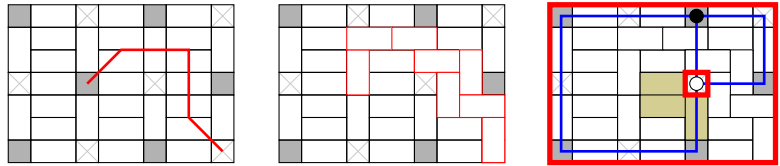
\begin{figure}
		\centering
		\begin{tikzpicture}[scale=0.3]
			\foreach \x in {-12,-6,...,18}
			{
				\foreach \y in {0,6} \draw [fill=gray!60] (\x,\y)--(\x,\y+1)--(\x+1,\y+1)--(\x+1,\y)--cycle;
			}
			\foreach \x in {-9,-3,...,21}
			{
				\foreach \y in {3} \draw [fill=gray!60] (\x,\y)--(\x,\y+1)--(\x+1,\y+1)--(\x+1,\y)--cycle;
			}
			\foreach \x in {-12,-6,...,18}
			{
				\foreach \y in {3}
				{
					\draw [color=gray!40] (\x,\y)--(\x+1,\y+1);
					\draw [color=gray!40] (\x,\y+1)--(\x+1,\y);
					\draw (\x,\y)--(\x,\y+1)--(\x+1,\y+1)--(\x+1,\y)--cycle;
				}
			}
			\foreach \x in {-9,-3,...,21}
			{
				\foreach \y in {0,6}
				{
					\draw [color=gray!40] (\x,\y)--(\x+1,\y+1);
					\draw [color=gray!40] (\x,\y+1)--(\x+1,\y);
					\draw (\x,\y)--(\x,\y+1)--(\x+1,\y+1)--(\x+1,\y)--cycle;
				}
			}
			\foreach \x in {-11,-8,-5,1,4,7,13,16,19}
			{
				\foreach \y in {0,...,6}
				{
					\draw (\x,\y)--(\x+2,\y)--(\x+2,\y+1)--(\x,\y+1)--cycle;
				}
			}
			\foreach \x in {-12,-9,...,21}
			{
				\foreach \y in {1,4}
				{
					\draw (\x,\y)--(\x+1,\y)--(\x+1,\y+2)--(\x,\y+2)--cycle;
				}
			}
			\draw [color=red,fill=white] (3,3)--(4,3)--(4,5)--(3,5)--cycle;
			\draw [color=red,fill=white] (3,5)--(5,5)--(5,6)--(3,6)--cycle;
			\draw [color=red,fill=white] (5,5)--(7,5)--(7,6)--(5,6)--cycle;
			\draw [color=red,fill=white] (6,4)--(6,5)--(8,5)--(8,4)--cycle;
			\draw [color=red,fill=white] (8,3)--(8,5)--(9,5)--(9,3)--cycle;
			\draw [color=red,fill=white] (7,2)--(7,4)--(8,4)--(8,2)--cycle;
			\draw [color=red,fill=white] (8,2)--(8,3)--(10,3)--(10,2)--cycle;
			\draw [color=red,fill=white] (9,0)--(9,2)--(10,2)--(10,0)--cycle;
			\draw [fill=white] (15,3)--(16,3)--(16,5)--(15,5)--cycle;
			\draw [fill=white] (15,5)--(17,5)--(17,6)--(15,6)--cycle;
			\draw [fill=white] (17,5)--(19,5)--(19,6)--(17,6)--cycle;
			\draw [fill=white] (18,4)--(18,5)--(20,5)--(20,4)--cycle;
			\draw [fill=white] (20,3)--(20,5)--(21,5)--(21,3)--cycle;
			\draw [fill=white] (19,2)--(19,4)--(20,4)--(20,2)--cycle;
			\draw [fill=white] (20,2)--(20,3)--(22,3)--(22,2)--cycle;
			\draw [fill=white] (21,0)--(21,2)--(22,2)--(22,0)--cycle;
			\draw [fill=olive!40] (18,2)--(18,4)--(16,4)--(16,2)--(18,2)--(18,1)--(19,1)--(19,3)--(16,3);
			\draw [color=blue,line width=1pt] (18.5,3.5)--(21.5,3.5)--(21.5,6.5)--(12.5,6.5)--(12.5,0.5)--(18.5,0.5)--(18.5,6.5);
			\draw [fill=white] (18.5,3.5) circle (0.3);
			\draw [fill=black] (18.5,6.5) circle (0.3);
			\draw [color=red,line width=2pt] (18,3)--(18,4)--(19,4)--(19,3)--cycle;
			\draw [color=red,line width=2pt] (12,0)--(22,0)--(22,7)--(12,7)--cycle;
			\draw [color=red,line width=1pt] (-8.5,3.5)--(-7,5)--(-4,5)--(-4,2)--(-2.5,0.5);
		\end{tikzpicture}
		\captionof{figure}{First graph from fig 10: expanded, substituted, with regions}
	\end{figure}

	Now, given a Hamiltonian cycle in the graph, we can clearly draw a loop that visits some each white monomino once, and the other polyominoes once or no times. Then, we can replace the used regions with the first, second, third and seventh polyominoes in figure 13. Then for polyominoes on the graph but not on the path, we can pair alternating white/black monominoes with the dominoes between them, and use the fourth and fifth polyominoes in figure 13, and for any remaining dominoes use the sixth polyomino. This allows us to construct a collection of loops passing through every cell, including one loop which represents the Hamiltonian cycle. Then observe that we can join these loops together by connecting them at adjoining blue edges (as long as they aren't across region boundaries), thus creating a single loop going through every cell once that visits each region once.
		
	\begin{figure}
		\centering
		\begin{tikzpicture}[scale=1/11]
			\draw [line width=1.5pt] (0,0)--(0,10)--(10,10)--(10,0)--cycle;
			\foreach \x in {2,4,6,8} \draw (\x,0)--(\x,10);
			\foreach \y in {2,4,6,8} \draw (0,\y)--(10,\y);
			\draw [color=red,line width=2pt] (5,-1)--(5,1)--(1,1)--(1,3)--(7,3)--(7,1)--(9,1)--(9,2);
			\draw [color=blue,line width=2pt] (9,2)--(9,8);
			\draw [color=red,line width=2pt] (9,8)--(9,9)--(8,9);
			\draw [color=blue,line width=2pt] (8,9)--(2,9);
			\draw [color=red,line width=2pt] (2,9)--(1,9)--(1,7)--(7,7)--(7,5)--(-1,5);
				
			\draw [line width=1.5pt] (12,0)--(12,10)--(22,10)--(22,0)--cycle;
			\foreach \x in {14,16,18,20} \draw (\x,0)--(\x,10);
			\foreach \y in {2,4,6,8} \draw (12,\y)--(22,\y);
			\draw [color=red,line width=2pt] (17,-1)--(17,1)--(13,1)--(13,2);
			\draw [color=blue,line width=2pt] (13,2)--(13,8);
			\draw [color=red,line width=2pt] (13,8)--(13,9)--(15,9)--(15,3)--(19,3)--(19,1)--(21,1)--(21,2);
			\draw [color=blue,line width=2pt] (21,2)--(21,8);
			\draw [color=red,line width=2pt] (21,8)--(21,9)--(19,9)--(19,5)--(17,5)--(17,11);
				
			\draw [line width=1.5pt] (24,0)--(24,20)--(34,20)--(34,0)--cycle;
			\foreach \x in {26,28,30,32} \draw (\x,0)--(\x,20);
			\foreach \y in {2,4,...,18} \draw (24,\y)--(34,\y);
			\draw [color=red,line width=2pt] (29,-1)--(29,1)--(25,1)--(25,2);
			\draw [color=blue,line width=2pt] (25,2)--(25,8);
			\draw [color=red,line width=2pt] (25,8)--(25,9)--(27,9)--(27,3)--(31,3)--(31,1)--(33,1)--(33,2);
			\draw [color=blue,line width=2pt] (33,2)--(33,8);
			\draw [color=red,line width=2pt] (33,8)--(33,9)--(31,9)--(31,5)--(29,5)--(29,11)--(25,11)--(25,12);
			\draw [color=blue,line width=2pt] (25,12)--(25,18);
			\draw [color=red,line width=2pt] (25,18)--(25,19)--(27,19)--(27,13)--(31,13)--(31,11)--(33,11)--(33,12);
			\draw [color=blue,line width=2pt] (33,12)--(33,18);
			\draw [color=red,line width=2pt] (33,18)--(33,19)--(31,19)--(31,15)--(29,15)--(29,21);
				
			\draw [line width=1.5pt] (36,30)--(46,30)--(46,0)--(36,0)--(36,40)--(46,40)--(46,10)--(36,10);
			\foreach \x in {38,40,42,44} \draw (\x,0)--(\x,40);
			\foreach \y in {2,4,...,38} \draw (36,\y)--(46,\y);
			\draw [color=red,line width=2pt] (38,1)--(37,1)--(37,2);
			\draw [color=blue,line width=2pt] (37,2)--(37,8);
			\draw [color=red,line width=2pt] (37,8)--(37,12);
			\draw [color=blue,line width=2pt] (37,12)--(37,18);
			\draw [color=red,line width=2pt] (37,18)--(37,19)--(39,19)--(39,3)--(43,3)--(43,5)--(41,5)--(41,7)--(41,7)--(43,7)--(43,9)--(41,9)--(41,11)--(43,11)--(43,13)--(41,13)--(41,15)--(43,15)--(43,17)--(41,17)--(41,19)--(43,19)--(43,21)--(41,21)--(41,23)--(43,23)--(43,25)--(41,25)--(41,27)--(43,27)--(43,29)--(41,29)--(41,31)--(43,31)--(43,33)--(41,33)--(41,35)--(43,35)--(43,37)--(39,37)--(39,21)--(37,21)--(37,22);
			\draw [color=blue,line width=2pt] (37,22)--(37,28);
			\draw [color=red,line width=2pt] (37,28)--(37,32);
			\draw [color=blue,line width=2pt] (37,32)--(37,38);
			\draw [color=red,line width=2pt] (37,38)--(37,39)--(38,39);
			\draw [color=blue,line width=2pt] (38,39)--(44,39);
			\draw [color=red,line width=2pt] (44,39)--(45,39)--(45,38);
			\draw [color=blue,line width=2pt] (45,38)--(45,32);
			\draw [color=red,line width=2pt] (45,32)--(45,28);
			\draw [color=blue,line width=2pt] (45,28)--(45,22);
			\draw [color=red,line width=2pt] (45,22)--(45,18);
			\draw [color=blue,line width=2pt] (45,18)--(45,12);
			\draw [color=red,line width=2pt] (45,12)--(45,8);
			\draw [color=blue,line width=2pt] (45,8)--(45,2);
			\draw [color=red,line width=2pt] (45,2)--(45,1)--(44,1);
			\draw [color=blue,line width=2pt] (44,1)--(38,1);
				
			\draw [line width=1.5pt] (68,20)--(68,0)--(58,0)--(58,10)--(48,10)--(48,20)--(78,20)--(78,30)--(68,30)--(68,40)--(58,40)--(58,20);
			\draw [line width=1.5pt] (58,30)--(68,30);
			\draw [line width=1.5pt] (58,10)--(68,10);
			\foreach \x in {60,62,64,66} \draw (\x,0)--(\x,40);
			\foreach \x in {50,52,...,58} \draw (\x,10)--(\x,20);
			\foreach \x in {68,70,...,76} \draw (\x,20)--(\x,30);
			\foreach \y in {2,4,6,8,32,34,36,38} \draw (58,\y)--(68,\y);
			\foreach \y in {12,14,16,18} \draw (48,\y)--(68,\y);
			\foreach \y in {22,24,26,28} \draw (58,\y)--(78,\y);
			\draw [color=red,line width=2pt] (67,12)--(67,8);
			\draw [color=red,line width=2pt] (67,2)--(67,1)--(66,1);
			\draw [color=red,line width=2pt] (60,1)--(59,1)--(59,2);
			\draw [color=red,line width=2pt] (59,8)--(59,15)--(57,15)--(57,17)--(55,17)--(55,15)--(53,15)--(53,17)--(51,17)--(51,13)--(57,13)--(57,11)--(56,11);
			\draw [color=red,line width=2pt] (50,11)--(49,11)--(49,12);
			\draw [color=red,line width=2pt] (49,18)--(49,19)--(50,19);
			\draw [color=red,line width=2pt] (56,19)--(59,19)--(59,17)--(61,17)--(61,3)--(65,3)--(65,5)--(63,5)--(63,7)--(65,7)--(65,9)--(63,9)--(63,11)--(65,11)--(65,13)--(63,13)--(63,15)--(65,15)--(65,17)--(63,17)--(63,19)--(61,19)--(61,21)--(59,21)--(59,22);
			\draw [color=red,line width=2pt] (59,28)--(59,32);
			\draw [color=red,line width=2pt] (59,38)--(59,39)--(60,39);
			\draw [color=red,line width=2pt] (66,39)--(67,39)--(67,38);
			\draw [color=red,line width=2pt] (67,32)--(67,25)--(69,25)--(69,23)--(71,23)--(71,25)--(73,25)--(73,23)--(75,23)--(75,27)--(69,27)--(69,29)--(70,29);
			\draw [color=red,line width=2pt] (76,29)--(77,29)--(77,28);
			\draw [color=red,line width=2pt] (77,22)--(77,21)--(76,21);
			\draw [color=red,line width=2pt] (70,21)--(67,21)--(67,23)--(65,23)--(65,37)--(61,37)--(61,35)--(63,35)--(63,33)--(61,33)--(61,31)--(63,31)--(63,29)--(61,29)--(61,27)--(63,27)--(63,25)--(61,25)--(61,23)--(63,23)--(63,21)--(65,21)--(65,19)--(67,19)--(67,18);
			\draw [color=blue,line width=2pt] (49,12)--(49,18);
			\draw [color=blue,line width=2pt] (59,2)--(59,8);
			\draw [color=blue,line width=2pt] (59,22)--(59,28);
			\draw [color=blue,line width=2pt] (59,32)--(59,38);
			\draw [color=blue,line width=2pt] (67,2)--(67,8);
			\draw [color=blue,line width=2pt] (67,12)--(67,18);
			\draw [color=blue,line width=2pt] (67,32)--(67,38);
			\draw [color=blue,line width=2pt] (77,22)--(77,28);
			\draw [color=blue,line width=2pt] (50,11)--(56,11);
			\draw [color=blue,line width=2pt] (50,19)--(56,19);
			\draw [color=blue,line width=2pt] (60,1)--(66,1);
			\draw [color=blue,line width=2pt] (60,39)--(66,39);
			\draw [color=blue,line width=2pt] (70,21)--(76,21);
			\draw [color=blue,line width=2pt] (70,29)--(76,29);
				
			\draw [line width=1.5pt] (80,0)--(80,20)--(90,20)--(90,0)--cycle;
			\foreach \x in {82,84,86,88} \draw (\x,0)--(\x,20);
			\foreach \y in {2,4,...,18} \draw (80,\y)--(90,\y);
			\draw [color=red,line width=2pt] (82,1)--(81,1)--(81,2);
			\draw [color=red,line width=2pt] (88,1)--(89,1)--(89,2);
			\draw [color=red,line width=2pt] (89,8)--(89,12);
			\draw [color=red,line width=2pt] (89,18)--(89,19)--(88,19);
			\draw [color=red,line width=2pt] (82,19)--(81,19)--(81,18);
			\draw [color=red,line width=2pt] (81,8)--(81,9)--(83,9)--(83,3)--(87,3)--(87,5)--(85,5)--(85,7)--(87,7)--(87,9)--(85,9)--(85,11)--(87,11)--(87,13)--(85,13)--(85,15)--(87,15)--(87,17)--(83,17)--(83,11)--(81,11)--(81,12);
			\draw [color=blue,line width=2pt] (81,2)--(81,8);
			\draw [color=blue,line width=2pt] (81,12)--(81,18);
			\draw [color=blue,line width=2pt] (89,2)--(89,8);
			\draw [color=blue,line width=2pt] (89,12)--(89,18);
			\draw [color=blue,line width=2pt] (82,1)--(88,1);
			\draw [color=blue,line width=2pt] (82,19)--(88,19);
				
			\draw [line width=1.5pt] (112,10)--(92,10)--(92,0)--(112,0)--(112,10)--(122,10)--(122,20)--(102,20)--(102,10);
			\foreach \x in {104,106,108,110} \draw (\x,0)--(\x,20);
			\foreach \x in {94,96,...,102} \draw (\x,0)--(\x,10);
			\foreach \x in {112,114,...,120} \draw (\x,10)--(\x,20);
			\foreach \y in {2,4,6,8,10} \draw (92,\y)--(112,\y);
			\foreach \y in {12,14,16,18} \draw (102,\y)--(122,\y);
			\draw [color=red,line width=2pt] (107,-1)--(107,1)--(100,1);
			\draw [color=red,line width=2pt] (94,1)--(93,1)--(93,2);
			\draw [color=red,line width=2pt] (93,8)--(93,9)--(94,9);
			\draw [color=red,line width=2pt] (100,9)--(103,9)--(103,7)--(95,7)--(95,3)--(97,3)--(97,5)--(99,5)--(99,3)--(101,3)--(101,5)--(103,5)--(103,3)--(105,3)--(105,7)--(109,7)--(109,5)--(107,5)--(107,3)--(109,3)--(109,1)--(111,1)--(111,2);
			\draw [color=red,line width=2pt] (111,8)--(111,9)--(109,9)--(109,11)--(107,11)--(107,9)--(105,9)--(105,11)--(103,11)--(103,12);
			\draw [color=red,line width=2pt] (103,18)--(103,19)--(105,19)--(105,17)--(107,17)--(107,15)--(105,15)--(105,13)--(109,13)--(109,17)--(111,17)--(111,15)--(113,15)--(113,17)--(115,17)--(115,15)--(117,15)--(117,17)--(119,17)--(119,13)--(111,13)--(111,11)--(114,11);
			\draw [color=red,line width=2pt] (120,11)--(121,11)--(121,12);
			\draw [color=red,line width=2pt] (121,18)--(121,19)--(120,19);
			\draw [color=red,line width=2pt] (114,19)--(107,19)--(107,21);
			\draw [color=blue,line width=2pt] (94,1)--(100,1);
			\draw [color=blue,line width=2pt] (94,9)--(100,9);
			\draw [color=blue,line width=2pt] (114,11)--(120,11);
			\draw [color=blue,line width=2pt] (114,19)--(120,19);
			\draw [color=blue,line width=2pt] (93,2)--(93,8);
			\draw [color=blue,line width=2pt] (121,12)--(121,18);
			\draw [color=blue,line width=2pt] (111,2)--(111,8);
			\draw [color=blue,line width=2pt] (103,12)--(103,18);
		\end{tikzpicture}
		\captionof{figure}{Possible paths for various polyomino combinations}
	\end{figure}
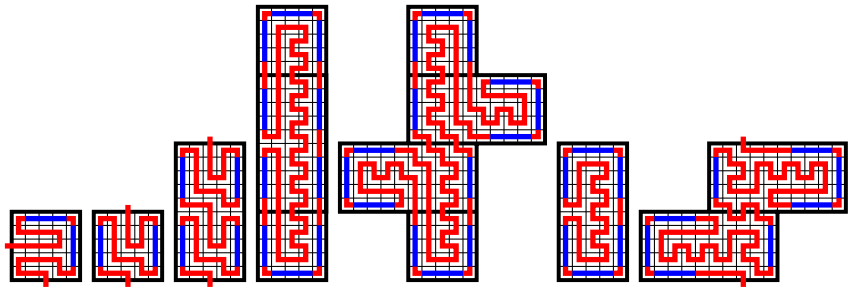
		
	To turn this into a path, we simply place a finish in a corner with a start next to it. To ensure each region is visited only once, we place `1' clues along the borders of each region except the one with the start and finish: in that region we place a `2' next to the finish (so the region is exited once then entered once). Finally, by observing that we can't visit two white regions in a row, (since due to their extra white cell from the unpaired vertex monomino, any visit must start and end on a white cell) so visiting two consecutive white regions would require visiting two white cells in a row which is impossible. So a path is only possible to make if alternates white and black regions, and thus only possible if a hamiltonian cycle exists in the graph. Thus, we are able to construct an appropriate Haisu puzzle, as shown in figure 14 (vastly simplified, though, since many intermediate steps have been omitted).
		
	\begin{figure}
		\centering
		\begin{tikzpicture}[scale=0.1,font=\sffamily]
			\foreach \x in {2,4,...,98} \draw (\x,0)--(\x,70);
			\foreach \y in {2,4,...,68} \draw (0,\y)--(100,\y);
			\draw [line width=2pt] (0,0)--(100,0)--(100,70)--(0,70)--cycle;
			\draw [line width=2pt] (60,30)--(60,40)--(70,40)--(70,30)--cycle;
			\foreach \x in {61,63,65,67,69} \draw (\x,31) node {\tiny 1};
			\foreach \x in {61,63,65,67,69} \draw (\x,38.9) node {\tiny 1};
			\foreach \y in {33,35,37} \draw (61,\y) node {\tiny 1};
			\foreach \y in {33,35,37} \draw (69,\y) node {\tiny 1};
			\draw [line width=2pt] (2,2)--(0,2);
			\draw [line width=2pt] (0,4)--(2,4)--(2,0);
			\draw (1,1) node {\scalebox{0.4}{F}};
			\draw (1,3) node {\scalebox{0.4}{S}};
			\draw (3,1) node {\tiny 2};
		\end{tikzpicture}
		\begin{tikzpicture}[scale=0.1,font=\sffamily]
		\foreach \x in {2,4,...,98} \draw (\x,0)--(\x,70);
		\foreach \y in {2,4,...,68} \draw (0,\y)--(100,\y);
		\draw [color=blue,line width=1pt] (10,0)--(10,70);
		\draw [color=blue,line width=1pt] (30,0)--(30,70);
		\draw [color=blue,line width=1pt] (40,0)--(40,50);
		\draw [color=blue,line width=1pt] (40,60)--(40,70);
		\draw [color=blue,line width=1pt] (50,50)--(50,60);
		\draw [color=blue,line width=1pt] (60,0)--(60,50);
		\draw [color=blue,line width=1pt] (60,60)--(60,70);
		\draw [color=blue,line width=1pt] (70,0)--(70,40);
		\draw [color=blue,line width=1pt] (70,50)--(70,70);
		\draw [color=blue,line width=1pt] (80,20)--(80,50);
		\draw [color=blue,line width=1pt] (90,0)--(90,20);
		\draw [color=blue,line width=1pt] (90,30)--(90,70);
		\draw [color=blue,line width=1pt] (0,10)--(90,10);
		\draw [color=blue,line width=1pt] (10,20)--(30,20);
		\draw [color=blue,line width=1pt] (40,20)--(60,20);
		\draw [color=blue,line width=1pt] (70,20)--(100,20);
		\draw [color=blue,line width=1pt] (0,30)--(70,30);
		\draw [color=blue,line width=1pt] (80,30)--(100,30);
		\draw [color=blue,line width=1pt] (0,40)--(30,40);
		\draw [color=blue,line width=1pt] (40,40)--(80,40);
		\draw [color=blue,line width=1pt] (90,40)--(100,40);
		\draw [color=blue,line width=1pt] (10,50)--(90,50);
		\draw [color=blue,line width=1pt] (0,60)--(100,60);
		\draw [line width=2pt] (0,0)--(100,0)--(100,70)--(0,70)--cycle;
		\draw [line width=2pt] (60,30)--(60,40)--(70,40)--(70,30)--cycle;
		\foreach \x in {61,63,65,67,69} \draw (\x,31) node {\tiny 1};
		\foreach \x in {61,63,65,67,69} \draw (\x,38.9) node {\tiny 1};
		\foreach \y in {33,35,37} \draw (61,\y) node {\tiny 1};
		\foreach \y in {33,35,37} \draw (69,\y) node {\tiny 1};
		\draw [line width=2pt] (2,2)--(0,2);
		\draw [line width=2pt] (0,4)--(2,4)--(2,0);
		\draw (1,1) node {\scalebox{0.4}{F}};
		\draw (1,3) node {\scalebox{0.4}{S}};
		\draw (3,1) node {\tiny 2};
		\draw [color=red,line width=2pt] (1,1)--(9,1)--(9,3)--(3,3)--(3,5)--(11,5)--(11,9)--(19,9)--(19,7)--(13,7)--(13,3)--(11,3)--(11,1)--(19,1)--(19,3)--(15,3)--(15,5)--(21,5)--(21,9)--(29,9)--(29,7)--(23,7)--(23,3)--(21,3)--(21,1)--(29,1)--(29,3)--(25,3)--(25,5)--(31,5)--(31,9)--(39,9)--(39,7)--(33,7)--(33,3)--(31,3)--(31,1)--(39,1)--(39,3)--(35,3)--(35,5)--(41,5)--(41,9)--(49,9)--(49,7)--(43,7)--(43,3)--(41,3)--(41,1)--(49,1)--(49,3)--(45,3)--(45,5)--(51,5)--(51,9)--(59,9)--(59,7)--(53,7)--(53,3)--(51,3)--(51,1)--(59,1)--(59,3)--(55,3)--(55,5)--(61,5)--(61,9)--(63,9)--(63,3)--(61,3)--(61,1)--(69,1)--(69,9)--(67,9)--(67,3)--(65,3)--(65,11)--(61,11)--(61,19)--(63,19)--(63,13)--(67,13)--(67,11)--(69,11)--(69,19)--(67,19)--(67,15)--(65,15)--(65,21)--(61,21)--(61,29)--(63,29)--(63,23)--(67,23)--(67,21)--(69,21)--(69,29)--(67,29)--(67,25)--(65,25)--(65,31)--(61,31)--(61,39)--(63,39)--(63,33)--(67,33)--(67,31)--(69,31)--(69,39)--(67,39)--(67,35)--(65,35)--(65,41)--(79,41)--(79,49)--(69,49)--(69,47)--(77,47)--(77,43)--(75,43)--(75,45)--(73,45)--(73,43)--(71,43)--(71,45)--(69,45)--(69,43)--(67,43)--(67,47)--(63,47)--(63,45)--(65,45)--(65,43)--(63,43)--(63,41)--(61,41)--(61,49)--(63,49)--(63,51)--(65,51)--(65,49)--(67,49)--(67,51)--(69,51)--(69,59)--(67,59)--(67,57)--(65,57)--(65,55)--(67,55)--(67,53)--(63,53)--(63,57)--(61,57)--(61,55)--(59,55)--(59,57)--(57,57)--(57,55)--(55,55)--(55,57)--(53,57)--(53,53)--(61,53)--(61,51)--(51,51)--(51,59)--(65,59)--(65,61)--(61,61)--(61,63)--(67,63)--(67,61)--(69,61)--(69,69)--(61,69)--(61,67)--(67,67)--(67,65)--(55,65)--(55,63)--(59,63)--(59,61)--(51,61)--(51,63)--(53,63)--(53,67)--(59,67)--(59,69)--(51,69)--(51,65)--(45,65)--(45,63)--(49,63)--(49,61)--(41,61)--(41,63)--(43,63)--(43,67)--(49,67)--(49,69)--(41,69)--(41,65)--(35,65)--(35,63)--(39,63)--(39,61)--(31,61)--(31,63)--(33,63)--(33,67)--(39,67)--(39,69)--(31,69)--(31,65)--(25,65)--(25,63)--(29,63)--(29,61)--(21,61)--(21,63)--(23,63)--(23,67)--(29,67)--(29,69)--(21,69)--(21,65)--(15,65)--(15,63)--(19,63)--(19,61)--(11,61)--(11,63)--(13,63)--(13,67)--(19,67)--(19,69)--(11,69)--(11,65)--(9,65)--(9,61)--(7,61)--(7,67)--(9,67)--(9,69)--(1,69)--(1,61)--(3,61)--(3,67)--(5,67)--(5,55)--(7,55)--(7,59)--(9,59)--(9,51)--(7,51)--(7,53)--(3,53)--(3,59)--(1,59)--(1,51)--(5,51)--(5,45)--(7,45)--(7,49)--(9,49)--(9,41)--(7,41)--(7,43)--(3,43)--(3,49)--(1,49)--(1,41)--(5,41)--(5,35)--(7,35)--(7,39)--(9,39)--(9,31)--(7,31)--(7,33)--(3,33)--(3,39)--(1,39)--(1,31)--(5,31)--(5,25)--(7,25)--(7,29)--(9,29)--(9,21)--(7,21)--(7,23)--(3,23)--(3,29)--(1,29)--(1,21)--(5,21)--(5,15)--(7,15)--(7,19)--(9,19)--(9,11)--(7,11)--(7,13)--(3,13)--(3,19)--(1,19)--(1,11)--(5,11)--(5,9)--(9,9)--(9,7)--(3,7)--(3,9)--(1,9)--(1,3);
		\foreach \y in {10,20,30,40,50} \draw [color=brown,line width=1pt] (19,\y+5)--(19,\y+3)--(17,\y+3)--(17,\y+5)--(15,\y+5)--(15,\y+3)--(13,\y+3)--(13,\y+7)--(19,\y+7)--(19,\y+9)--(11,\y+9)--(11,\y+1)--(29,\y+1)--(29,\y+9)--(21,\y+9)--(21,\y+7)--(27,\y+7)--(27,\y+3)--(25,\y+3)--(25,\y+5)--(23,\y+5)--(23,\y+3)--(21,\y+3)--(21,\y+5)--cycle;
		\foreach \y in {10,20,30,40} \draw [color=brown,line width=1pt] (49,\y+5)--(49,\y+3)--(47,\y+3)--(47,\y+5)--(45,\y+5)--(45,\y+3)--(43,\y+3)--(43,\y+7)--(49,\y+7)--(49,\y+9)--(41,\y+9)--(41,\y+1)--(59,\y+1)--(59,\y+9)--(51,\y+9)--(51,\y+7)--(57,\y+7)--(57,\y+3)--(55,\y+3)--(55,\y+5)--(53,\y+5)--(53,\y+3)--(51,\y+3)--(51,\y+5)--cycle;
		\draw [color=brown,line width=1pt] (39,55)--(39,53)--(37,53)--(37,55)--(35,55)--(35,53)--(33,53)--(33,57)--(39,57)--(39,59)--(31,59)--(31,51)--(49,51)--(49,59)--(41,59)--(41,57)--(47,57)--(47,53)--(45,53)--(45,55)--(43,55)--(43,53)--(41,53)--(41,55)--cycle;
		\foreach \y in {0,10,50,60} \draw [color=brown,line width=1pt] (79,\y+5)--(79,\y+3)--(77,\y+3)--(77,\y+5)--(75,\y+5)--(75,\y+3)--(73,\y+3)--(73,\y+7)--(79,\y+7)--(79,\y+9)--(71,\y+9)--(71,\y+1)--(89,\y+1)--(89,\y+9)--(81,\y+9)--(81,\y+7)--(87,\y+7)--(87,\y+3)--(85,\y+3)--(85,\y+5)--(83,\y+5)--(83,\y+3)--(81,\y+3)--(81,\y+5)--cycle;
		\draw [color=brown,line width=1pt] (89,25)--(89,23)--(87,23)--(87,25)--(85,25)--(85,23)--(83,23)--(83,27)--(89,27)--(89,29)--(81,29)--(81,21)--(99,21)--(99,29)--(91,29)--(91,27)--(97,27)--(97,23)--(95,23)--(95,25)--(93,25)--(93,23)--(91,23)--(91,25)--cycle;
		\draw [color=brown,line width=1pt] (31,11)--(39,11)--(39,29)--(31,29)--(31,21)--(33,21)--(33,27)--(37,27)--(37,25)--(35,25)--(35,23)--(37,23)--(37,21)--(35,21)--(35,19)--(37,19)--(37,17)--(35,17)--(35,15)--(37,15)--(37,13)--(33,13)--(33,19)--(31,19)--cycle;
		\draw [color=brown,line width=1pt] (31,31)--(39,31)--(39,49)--(31,49)--(31,41)--(33,41)--(33,47)--(37,47)--(37,45)--(35,45)--(35,43)--(37,43)--(37,41)--(35,41)--(35,39)--(37,39)--(37,37)--(35,37)--(35,35)--(37,35)--(37,33)--(33,33)--(33,39)--(31,39)--cycle;
		\draw [color=brown,line width=1pt] (81,31)--(89,31)--(89,49)--(81,49)--(81,41)--(83,41)--(83,47)--(87,47)--(87,45)--(85,45)--(85,43)--(87,43)--(87,41)--(85,41)--(85,39)--(87,39)--(87,37)--(85,37)--(85,35)--(87,35)--(87,33)--(83,33)--(83,39)--(81,39)--cycle;
		\draw [color=brown,line width=1pt] (71,21)--(79,21)--(79,39)--(71,39)--(71,31)--(73,31)--(73,37)--(77,37)--(77,35)--(75,35)--(75,33)--(77,33)--(77,31)--(75,31)--(75,29)--(77,29)--(77,27)--(75,27)--(75,25)--(77,25)--(77,23)--(73,23)--(73,29)--(71,29)--cycle;
		\draw [color=brown,line width=1pt] (91,1)--(99,1)--(99,19)--(91,19)--(91,11)--(93,11)--(93,17)--(97,17)--(97,15)--(95,15)--(95,13)--(97,13)--(97,11)--(95,11)--(95,9)--(97,9)--(97,7)--(95,7)--(95,5)--(97,5)--(97,3)--(93,3)--(93,9)--(91,9)--cycle;
		\draw [color=brown,line width=1pt] (91,31)--(99,31)--(99,69)--(91,69)--(91,51)--(93,51)--(93,67)--(97,67)--(97,65)--(95,65)--(95,63)--(97,63)--(97,61)--(95,61)--(95,59)--(97,59)--(97,57)--(95,57)--(95,55)--(97,55)--(97,53)--(95,53)--(95,51)--(97,51)--(97,49)--(95,49)--(95,47)--(97,47)--(97,45)--(95,45)--(95,43)--(97,43)--(97,41)--(95,41)--(95,39)--(97,39)--(97,37)--(95,37)--(95,35)--(97,35)--(97,33)--(93,33)--(93,49)--(91,49)--cycle;
		\foreach \y in {20,30,40,50,60} {\vtwist{14}{\y};}
		\foreach \y in {30,50,60} {\vtwist{34}{\y};}
		\foreach \y in {20,30,40,50} {\vtwist{54}{\y};}
		\foreach \y in {4,14,24,54,64} {\htwist{70}{\y};}
		\htwist{80}{24};
		\htwist{80}{34};
		\vtwist{94}{20};
		\vtwist{94}{30};
		\end{tikzpicture}
		\captionof{figure}{Example Haisu construction and solution using fig 12 (simplified)}
	\end{figure}
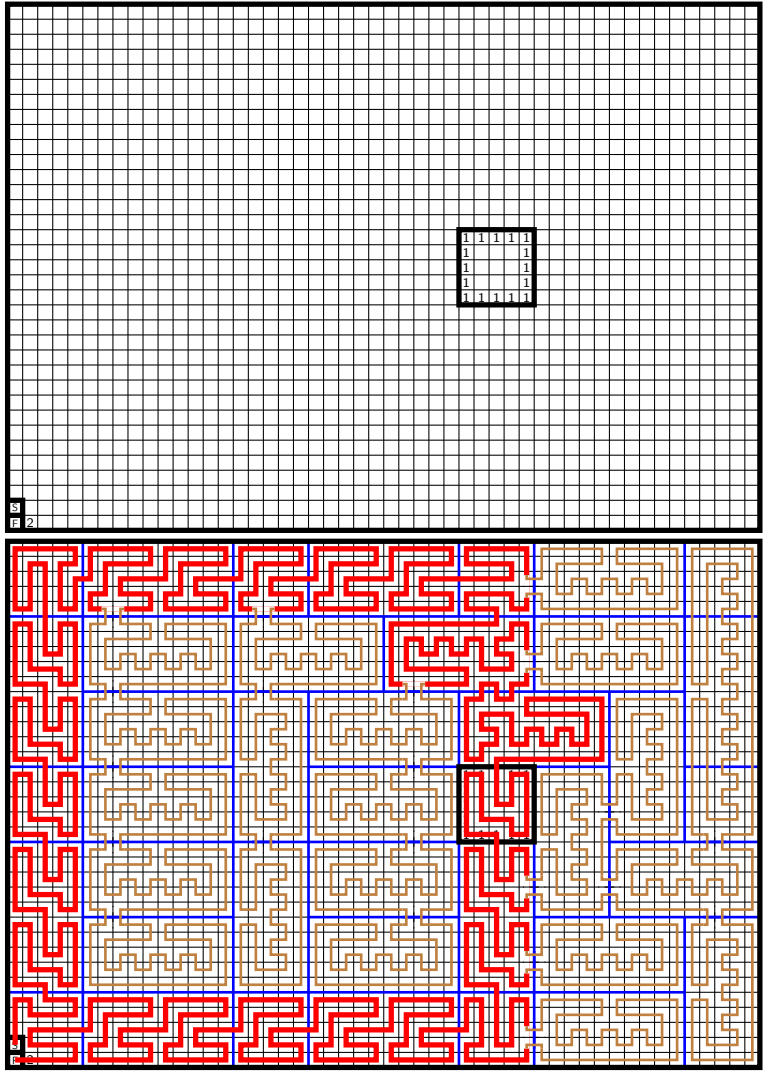
		
	To prove this is a polynomial mapping, observe that the only operations we did to the initial rectilinear realisation (which was $O(n^2)$ size) were dilating it by a constant factor several times, and unwrapping it which would have turned the graph into $O(n^4)$ size. Thus Haisu satisfiability is NP-complete.
		
	\subsection{Barred Simple Loop: Above Haisu reduction}
		
	To help with proving NP-completeness of the last two genres, we will introduce an intermediate genre: Barred Simple Loop. In Barred Simple Loop, the aim is to draw a loop going through every square, with the restriction that some edges between squares are marked and cannot be passed over.
		
	We are only going to reduce from the particular reduction previously exhibited. Recall that we wanted each region to be entered and exited exactly once, and for every region it was to be adjacent to exactly three regions of the opposite colour (the multigraph in figure 10 meant that the construction in figure 14 had an adjacency with 'multiplicity' 3, but that's irrelevant to our argument).
		
	As a result, we can simply put borders along the edges of all regions except for three openings in the middles of the sides of black regions to their appropriate white regions. This is because each region must be entered and exited the same number of times, thus the total number of entries/exits must be even, and since each region must be entered but there are only three points to do so, each region must be entered and exited exactly once. We can use the same constructions from Haisu to build a valid solution.
		
	This actually has the same size as the Haisu reduction, so is found in polynomial space. Furthermore it is clearly easy to check satisfaction in $O(n^2)$ time, thus Barred Simple Loop satisfiability must be NP-complete.
		
	\begin{figure}
		\centering
		\begin{tikzpicture}[scale=0.1]
			\foreach \x in {2,4,...,98} \draw (\x,0)--(\x,100);
			\foreach \y in {2,4,...,98} \draw (0,\y)--(100,\y);
			\foreach \x in {0,30,60,90}{\foreach \y in {0,30,60,90} \draw [color=blue,line width=1pt] (\x,\y)--(\x,\y+10)--(\x+10,\y+10)--(\x+10,\y)--cycle;}
			\foreach \x in {10,40,70}{\foreach \y in {0,10,...,90}\draw [color=blue,line width=1pt] (\x,\y)--(\x+20,\y)--(\x+20,\y+10)--(\x,\y+10)--cycle;}
			\draw [line width=2pt] (0,0)--(0,100)--(100,100)--(100,0)--cycle;
			\draw [line width=2pt] (60,100)--(60,96);
			\draw [line width=2pt] (60,94)--(60,90)--(64,90);
			\draw [line width=2pt] (66,90)--(70,90)--(70,94);
			\draw [line width=2pt] (70,100)--(70,96);
			\draw [line width=2pt] (34,60)--(30,60)--(30,70)--(34,70);
			\draw [line width=2pt] (36,70)--(40,70)--(40,66);
			\draw [line width=2pt] (40,64)--(40,60)--(36,60);
			\draw [line width=2pt] (66,40)--(70,40)--(70,30)--(66,30);
			\draw [line width=2pt] (64,30)--(60,30)--(60,34);
			\draw [line width=2pt] (60,36)--(60,40)--(64,40);
			\draw [line width=2pt] (30,0)--(30,4);
			\draw [line width=2pt] (30,6)--(30,10)--(34,10);
			\draw [line width=2pt] (36,10)--(40,10)--(40,6);
			\draw [line width=2pt] (40,4)--(40,0);
			\draw [line width=2pt] (30,10)--(30,60);
			\draw [line width=2pt] (70,40)--(70,90);
			\draw [line width=2pt] (40,10)--(40,20)--(60,20)--(60,30);
			\draw [line width=2pt] (60,40)--(60,50)--(40,50)--(40,60);
			\draw [line width=2pt] (40,70)--(40,80)--(60,80)--(60,90);
		\end{tikzpicture}
		\captionof{figure}{Example Barred Simple Loop construction from fig 8 graph}
	\end{figure}
	
	\subsection{Oriental House: Barred Simple Loop}
		
	Firstly, we notice that the left of figure 16 displays a region with a single solution up to the orientation of the loop. This is because if we consider the middle arrow, we can't go horizontally, thus we must pass through that cell and one of the other arrows. This means that we must pass through two arrows in opposite N--S directions. But the only way to satisfy this is to enter and leave the region both on north edges or both on south edges, and out of these only the north edge option works. This can also be extended as seen in that figure.
		
	\begin{figure}
		\centering
		\begin{tikzpicture}[scale=0.5]
			\draw [line width=2pt] (0,0)--(0,1)--(1,1)--(1,3)--(2,3)--(2,0)--cycle;
			\draw [line width=2pt] (3,0)--(3,1)--(4,1)--(4,3)--(5,3)--(5,0)--cycle;
			\draw [line width=2pt] (6,0)--(6,1)--(7,1)--(7,4)--(8,4)--(8,0)--cycle;
			\draw [line width=2pt] (9,0)--(9,1)--(10,1)--(10,4)--(11,4)--(11,0)--cycle;
			\draw [line width=2pt] (12,0)--(12,1)--(13,1)--(13,5)--(14,5)--(14,0)--cycle;
			\draw (1,0)--(1,1)--(2,1);
			\draw (1,2)--(2,2);
			\draw (4,0)--(4,1)--(5,1);
			\draw (4,2)--(5,2);
			\draw (7,0)--(7,1)--(8,1);
			\draw (7,2)--(8,2);
			\draw (7,3)--(8,3);
			\draw (10,0)--(10,1)--(11,1);
			\draw (10,2)--(11,2);
			\draw (10,3)--(11,3);
			\draw (13,0)--(13,1)--(14,1);
			\draw (13,2)--(14,2);
			\draw (13,3)--(14,3);
			\draw (13,4)--(14,4);
			\draw (15,2.5) node {$\cdots$};
			\draw (1.5,0.5) node {$\uparrow$};
			\draw (1.5,1.5) node {$\downarrow$};
			\draw (1.5,2.5) node {$\uparrow$};
			\draw (4.5,0.5) node {$\uparrow$};
			\draw (4.5,1.5) node {$\downarrow$};
			\draw (4.5,2.5) node {$\uparrow$};
			\draw (7.5,1.5) node {$\uparrow$};
			\draw (7.5,2.5) node {$\downarrow$};
			\draw (7.5,3.5) node {$\uparrow$};
			\draw (10.5,1.5) node {$\uparrow$};
			\draw (10.5,2.5) node {$\downarrow$};
			\draw (10.5,3.5) node {$\uparrow$};
			\draw (13.5,2.5) node {$\uparrow$};
			\draw (13.5,3.5) node {$\downarrow$};
			\draw (13.5,4.5) node {$\uparrow$};
			\draw [color=red,line width=2pt] (3.5,1.5)--(3.5,0.5)--(4.5,0.5)--(4.5,3.5);
			\draw [color=red,line width=2pt] (9.5,1.5)--(9.5,0.5)--(10.5,0.5)--(10.5,4.5);
		\end{tikzpicture}
		\captionof{figure}{Oriental House regions with unique unoriented solutions}
	\end{figure}
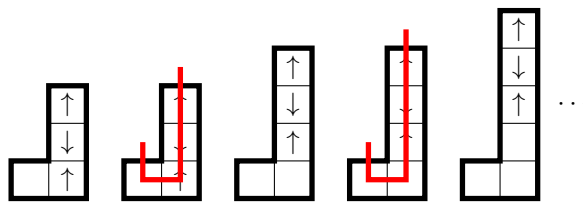
	
	Using this figure, we can construct a $13\times13$ `metacell' which acts like a normal cell except that we can block off any edges by adding the blue regions, pictured in figure 17. To show that this works, in figure 18, notice that the red lines are forced, and by adding the blue L-regions we can force the violet lines which block off the use of that edge. The left diagram illustrates a straight cell and the right a turning cell: due to the edgework we can only enter cells in the middles of the sides, and we can't have multiple entries into the one cell because, if we colour in a checkerboard so the entries are black, each visit must start and end on a black cell so the total number of used black cells minus the total number of used white cells is equal to the number of visits, and this is 1 because we must pass through every cell in the region.
	
	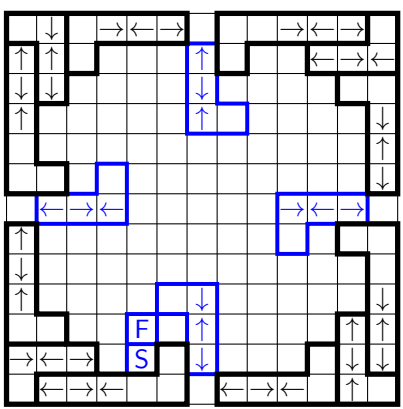
\begin{figure}
		\centering
		\begin{tikzpicture}[scale=0.2,font=\sffamily]
			\foreach \x in {0,2,...,26} \draw (\x,0)--(\x,26);
			\foreach \y in {0,2,...,26} \draw (0,\y)--(26,\y);
			\draw [color=blue,line width=1.5pt] (2,14)--(2,12)--(8,12)--(8,16)--(6,16)--(6,14)--(4,14);
			\draw [color=blue,line width=1.5pt] (8,2)--(8,6)--(10,6)--(10,4)--(8,4);
			\draw [color=blue,line width=1.5pt] (24,12)--(24,14)--(18,14)--(18,10)--(20,10)--(20,12)--(22,12);
			\draw [color=blue,line width=1.5pt] (12,2)--(14,2)--(14,8)--(10,8)--(10,6)--(12,6)--(12,4);
			\draw [color=blue,line width=1.5pt] (14,24)--(12,24)--(12,18)--(16,18)--(16,20)--(14,20)--(14,22);
			\draw [line width=2pt] (2,0)--(2,2)--(10,2)--(10,4)--(12,4)--(12,0)--(0,0)--(0,12)--(2,12)--(2,6)--(4,6)--(4,4)--(6,4)--(6,2);
			\draw [line width=2pt] (0,4)--(4,4);
			\draw [line width=2pt] (0,24)--(2,24)--(2,16)--(4,16)--(4,14)--(0,14)--(0,26)--(12,26)--(12,24)--(6,24)--(6,22)--(4,22)--(4,20)--(2,20);
			\draw [line width=2pt] (4,26)--(4,22);
			\draw [line width=2pt] (24,26)--(24,24)--(16,24)--(16,22)--(14,22)--(14,26)--(26,26)--(26,14)--(24,14)--(24,20)--(22,20)--(22,22)--(20,22)--(20,24);
			\draw [line width=2pt] (26,22)--(22,22);
			\draw [line width=2pt] (26,2)--(24,2)--(24,10)--(22,10)--(22,12)--(26,12)--(26,0)--(14,0)--(14,2)--(20,2)--(20,4)--(22,4)--(22,6)--(24,6);
			\draw [line width=2pt] (22,0)--(22,4);
			\draw (19,0.85) node {$\leftarrow$};
			\draw (17,0.85) node {$\rightarrow$};
			\draw (15,0.85) node {$\leftarrow$};
			\draw (7,0.85) node {$\leftarrow$};
			\draw (5,0.85) node {$\rightarrow$};
			\draw (3,0.85) node {$\leftarrow$};
			\draw (5,2.85) node {$\rightarrow$};
			\draw (3,2.85) node {$\leftarrow$};
			\draw (1,2.85) node {$\rightarrow$};
			\draw (1,7) node {$\uparrow$};
			\draw (1,9) node {$\downarrow$};
			\draw (1,11) node {$\uparrow$};
			\draw (1,19) node {$\uparrow$};
			\draw (1,21) node {$\downarrow$};
			\draw (1,23) node {$\uparrow$};
			\draw (3,21) node {$\downarrow$};
			\draw (3,23) node {$\uparrow$};
			\draw (3,25) node {$\downarrow$};
			\draw (7,24.85) node {$\rightarrow$};
			\draw (9,24.85) node {$\leftarrow$};
			\draw (11,24.85) node {$\rightarrow$};
			\draw (19,24.85) node {$\rightarrow$};
			\draw (21,24.85) node {$\leftarrow$};
			\draw (23,24.85) node {$\rightarrow$};
			\draw (21,22.85) node {$\leftarrow$};
			\draw (23,22.85) node {$\rightarrow$};
			\draw (25,22.85) node {$\leftarrow$};
			\draw (25,19) node {$\downarrow$};
			\draw (25,17) node {$\uparrow$};
			\draw (25,15) node {$\downarrow$};
			\draw (25,7) node {$\downarrow$};
			\draw (25,5) node {$\uparrow$};
			\draw (25,3) node {$\downarrow$};
			\draw (23,5) node {$\uparrow$};
			\draw (23,3) node {$\downarrow$};
			\draw (23,1) node {$\uparrow$};
			\draw [color=blue] (9,3) node {S};
			\draw [color=blue] (9,5) node {F};
			\draw [color=blue] (7,12.85) node {$\leftarrow$};
			\draw [color=blue] (5,12.85) node {$\rightarrow$};
			\draw [color=blue] (3,12.85) node {$\leftarrow$};
			\draw [color=blue] (23,12.85) node {$\rightarrow$};
			\draw [color=blue] (21,12.85) node {$\leftarrow$};
			\draw [color=blue] (19,12.85) node {$\rightarrow$};
			\draw [color=blue] (13,7) node {$\downarrow$};
			\draw [color=blue] (13,5) node {$\uparrow$};
			\draw [color=blue] (13,3) node {$\downarrow$};
			\draw [color=blue] (13,19) node {$\uparrow$};
			\draw [color=blue] (13,21) node {$\downarrow$};
			\draw [color=blue] (13,23) node {$\uparrow$};
		\end{tikzpicture}
		\captionof{figure}{Oriental House metacell}
	\end{figure}
	
	Finally, we can add a start/finish as indicated in any metacell, disconnecting the appropriate loop segment, without affecting the logic. Then, by using these metacells with appropriate blocks to form a Barred Simple Loop, we have a constant factor dilation of the Barred Simple Loop (which must therefore be polynomial) so Oriental House satisfiability is also NP-complete.
	
	\begin{figure}
		\centering
		\begin{tikzpicture}[scale=0.2,font=\sffamily]
			\foreach \x in {0,2,...,26} \draw (\x,0)--(\x,26);
			\foreach \y in {0,2,...,26} \draw (0,\y)--(26,\y);
			\draw [color=blue,line width=1.5pt] (2,14)--(2,12)--(8,12)--(8,16)--(6,16)--(6,14)--(4,14);
			\draw [color=blue,line width=1.5pt] (8,2)--(8,6)--(10,6)--(10,4)--(8,4);
			\draw [color=blue,line width=1.5pt] (24,12)--(24,14)--(18,14)--(18,10)--(20,10)--(20,12)--(22,12);
			\draw [line width=2pt] (2,0)--(2,2)--(10,2)--(10,4)--(12,4)--(12,0)--(0,0)--(0,12)--(2,12)--(2,6)--(4,6)--(4,4)--(6,4)--(6,2);
			\draw [line width=2pt] (0,4)--(4,4);
			\draw [line width=2pt] (0,24)--(2,24)--(2,16)--(4,16)--(4,14)--(0,14)--(0,26)--(12,26)--(12,24)--(6,24)--(6,22)--(4,22)--(4,20)--(2,20);
			\draw [line width=2pt] (4,26)--(4,22);
			\draw [line width=2pt] (24,26)--(24,24)--(16,24)--(16,22)--(14,22)--(14,26)--(26,26)--(26,14)--(24,14)--(24,20)--(22,20)--(22,22)--(20,22)--(20,24);
			\draw [line width=2pt] (26,22)--(22,22);
			\draw [line width=2pt] (26,2)--(24,2)--(24,10)--(22,10)--(22,12)--(26,12)--(26,0)--(14,0)--(14,2)--(20,2)--(20,4)--(22,4)--(22,6)--(24,6);
			\draw [line width=2pt] (22,0)--(22,4);
			\draw (19,0.85) node {$\leftarrow$};
			\draw (17,0.85) node {$\rightarrow$};
			\draw (15,0.85) node {$\leftarrow$};
			\draw (7,0.85) node {$\leftarrow$};
			\draw (5,0.85) node {$\rightarrow$};
			\draw (3,0.85) node {$\leftarrow$};
			\draw (5,2.85) node {$\rightarrow$};
			\draw (3,2.85) node {$\leftarrow$};
			\draw (1,2.85) node {$\rightarrow$};
			\draw (1,7) node {$\uparrow$};
			\draw (1,9) node {$\downarrow$};
			\draw (1,11) node {$\uparrow$};
			\draw (1,19) node {$\uparrow$};
			\draw (1,21) node {$\downarrow$};
			\draw (1,23) node {$\uparrow$};
			\draw (3,21) node {$\downarrow$};
			\draw (3,23) node {$\uparrow$};
			\draw (3,25) node {$\downarrow$};
			\draw (7,24.85) node {$\rightarrow$};
			\draw (9,24.85) node {$\leftarrow$};
			\draw (11,24.85) node {$\rightarrow$};
			\draw (19,24.85) node {$\rightarrow$};
			\draw (21,24.85) node {$\leftarrow$};
			\draw (23,24.85) node {$\rightarrow$};
			\draw (21,22.85) node {$\leftarrow$};
			\draw (23,22.85) node {$\rightarrow$};
			\draw (25,22.85) node {$\leftarrow$};
			\draw (25,19) node {$\downarrow$};
			\draw (25,17) node {$\uparrow$};
			\draw (25,15) node {$\downarrow$};
			\draw (25,7) node {$\downarrow$};
			\draw (25,5) node {$\uparrow$};
			\draw (25,3) node {$\downarrow$};
			\draw (23,5) node {$\uparrow$};
			\draw (23,3) node {$\downarrow$};
			\draw (23,1) node {$\uparrow$};
			\draw [color=blue] (9,3) node {S};
			\draw [color=blue] (9,5) node {F};
			\draw [color=blue] (7,12.85) node {$\leftarrow$};
			\draw [color=blue] (5,12.85) node {$\rightarrow$};
			\draw [color=blue] (3,12.85) node {$\leftarrow$};
			\draw [color=blue] (23,12.85) node {$\rightarrow$};
			\draw [color=blue] (21,12.85) node {$\leftarrow$};
			\draw [color=blue] (19,12.85) node {$\rightarrow$};
			\draw [color=red,line width=2pt] (9,4.5)--(9,3)--(11,3)--(11,1)--(1,1)--(1,3)--(7,3)--(7,5)--(5,5)--(5,6.5);
			\draw [color=red,line width=2pt] (3,6.5)--(3,5)--(1,5)--(1,12.5);
			\draw [color=red,line width=2pt] (4.5,17)--(3,17)--(3,15)--(1,15)--(1,25)--(3,25)--(3,19)--(5,19)--(5,21)--(6.5,21);
			\draw [color=red,line width=2pt] (6.5,23)--(5,23)--(5,25)--(12.5,25);
			\draw [color=red,line width=2pt] (17,21.5)--(17,23)--(15,23)--(15,25)--(25,25)--(25,23)--(19,23)--(19,21)--(21,21)--(21,19.5);
			\draw [color=red,line width=2pt] (23,19.5)--(23,21)--(25,21)--(25,13.5);
			\draw [color=red,line width=2pt] (21.5,9)--(23,9)--(23,11)--(25,11)--(25,1)--(23,1)--(23,7)--(21,7)--(21,5)--(19.5,5);
			\draw [color=red,line width=2pt] (19.5,3)--(21,3)--(21,1)--(13.5,1);
			\draw [color=violet,line width=1.5pt] (1,12.5)--(1,13)--(7,13)--(7,15)--(5,15)--(5,17)--(4.5,17);
			\draw [color=violet,line width=1.5pt] (25,13.5)--(25,13)--(19,13)--(19,11)--(21,11)--(21,9)--(21.5,9);
			\draw [color=brown,line width=1pt] (13,-1)--(13,1)--(13.5,1);
			\draw [color=brown,line width=1pt] (19.5,3)--(19,3)--(19,5)--(19.5,5);
			\draw [color=brown,line width=1pt] (23,19.5)--(23,19)--(21,19)--(21,19.5);
			\draw [color=brown,line width=1pt] (17,21.5)--(17,21)--(13,21)--(13,23)--(9,23)--(9,21)--(11,21)--(11,19)--(7,19)--(7,17)--(13,17)--(13,19)--(19,19)--(19,17)--(23,17)--(23,15)--(17,15)--(17,17)--(15,17)--(15,15)--(9,15)--(9,13)--(17,13)--(17,11)--(11,11)--(11,9)--(9,9)--(9,11)--(3,11)--(3,9)--(7,9)--(7,7)--(13,7)--(13,9)--(19,9)--(19,7)--(15,7)--(15,5)--(17,5)--(17,3)--(13,3)--(13,5)--(9,5)--(9,4.5);
			\draw [color=brown,line width=1pt] (5,6.5)--(5,7)--(3,7)--(3,6.5);
			\draw [color=brown,line width=1pt] (6.5,21)--(7,21)--(7,23)--(6.5,23);
			\draw [color=brown,line width=1pt] (12.5,25)--(13,25)--(13,27);
			
			\foreach \x in {0,2,...,26} \draw (\x+28,0)--(\x+28,26);
			\foreach \y in {0,2,...,26} \draw (0+28,\y)--(26+28,\y);
			\draw [color=blue,line width=1.5pt] (2+28,14)--(2+28,12)--(8+28,12)--(8+28,16)--(6+28,16)--(6+28,14)--(4+28,14);
			\draw [color=blue,line width=1.5pt] (8+28,2)--(8+28,6)--(10+28,6)--(10+28,4)--(8+28,4);
			\draw [color=blue,line width=1.5pt] (42,24)--(40,24)--(40,18)--(44,18)--(44,20)--(42,20)--(42,22);
			\draw [line width=2pt] (2+28,0)--(2+28,2)--(10+28,2)--(10+28,4)--(12+28,4)--(12+28,0)--(0+28,0)--(0+28,12)--(2+28,12)--(2+28,6)--(4+28,6)--(4+28,4)--(6+28,4)--(6+28,2);
			\draw [line width=2pt] (0+28,4)--(4+28,4);
			\draw [line width=2pt] (0+28,24)--(2+28,24)--(2+28,16)--(4+28,16)--(4+28,14)--(0+28,14)--(0+28,26)--(12+28,26)--(12+28,24)--(6+28,24)--(6+28,22)--(4+28,22)--(4+28,20)--(2+28,20);
			\draw [line width=2pt] (4+28,26)--(4+28,22);
			\draw [line width=2pt] (24+28,26)--(24+28,24)--(16+28,24)--(16+28,22)--(14+28,22)--(14+28,26)--(26+28,26)--(26+28,14)--(24+28,14)--(24+28,20)--(22+28,20)--(22+28,22)--(20+28,22)--(20+28,24);
			\draw [line width=2pt] (26+28,22)--(22+28,22);
			\draw [line width=2pt] (26+28,2)--(24+28,2)--(24+28,10)--(22+28,10)--(22+28,12)--(26+28,12)--(26+28,0)--(14+28,0)--(14+28,2)--(20+28,2)--(20+28,4)--(22+28,4)--(22+28,6)--(24+28,6);
			\draw [line width=2pt] (22+28,0)--(22+28,4);
			\draw (19+28,0.85) node {$\leftarrow$};
			\draw (17+28,0.85) node {$\rightarrow$};
			\draw (15+28,0.85) node {$\leftarrow$};
			\draw (7+28,0.85) node {$\leftarrow$};
			\draw (5+28,0.85) node {$\rightarrow$};
			\draw (3+28,0.85) node {$\leftarrow$};
			\draw (5+28,2.85) node {$\rightarrow$};
			\draw (3+28,2.85) node {$\leftarrow$};
			\draw (1+28,2.85) node {$\rightarrow$};
			\draw (1+28,7) node {$\uparrow$};
			\draw (1+28,9) node {$\downarrow$};
			\draw (1+28,11) node {$\uparrow$};
			\draw (1+28,19) node {$\uparrow$};
			\draw (1+28,21) node {$\downarrow$};
			\draw (1+28,23) node {$\uparrow$};
			\draw (3+28,21) node {$\downarrow$};
			\draw (3+28,23) node {$\uparrow$};
			\draw (3+28,25) node {$\downarrow$};
			\draw (7+28,24.85) node {$\rightarrow$};
			\draw (9+28,24.85) node {$\leftarrow$};
			\draw (11+28,24.85) node {$\rightarrow$};
			\draw (19+28,24.85) node {$\rightarrow$};
			\draw (21+28,24.85) node {$\leftarrow$};
			\draw (23+28,24.85) node {$\rightarrow$};
			\draw (21+28,22.85) node {$\leftarrow$};
			\draw (23+28,22.85) node {$\rightarrow$};
			\draw (25+28,22.85) node {$\leftarrow$};
			\draw (25+28,19) node {$\downarrow$};
			\draw (25+28,17) node {$\uparrow$};
			\draw (25+28,15) node {$\downarrow$};
			\draw (25+28,7) node {$\downarrow$};
			\draw (25+28,5) node {$\uparrow$};
			\draw (25+28,3) node {$\downarrow$};
			\draw (23+28,5) node {$\uparrow$};
			\draw (23+28,3) node {$\downarrow$};
			\draw (23+28,1) node {$\uparrow$};
			\draw [color=blue] (9+28,3) node {S};
			\draw [color=blue] (9+28,5) node {F};
			\draw [color=blue] (7+28,12.85) node {$\leftarrow$};
			\draw [color=blue] (5+28,12.85) node {$\rightarrow$};
			\draw [color=blue] (3+28,12.85) node {$\leftarrow$};
			\draw [color=blue] (41,23) node {$\uparrow$};
			\draw [color=blue] (41,21) node {$\downarrow$};
			\draw [color=blue] (41,19) node {$\uparrow$};
			\draw [color=red,line width=2pt] (9+28,4.5)--(9+28,3)--(11+28,3)--(11+28,1)--(1+28,1)--(1+28,3)--(7+28,3)--(7+28,5)--(5+28,5)--(5+28,6.5);
			\draw [color=red,line width=2pt] (3+28,6.5)--(3+28,5)--(1+28,5)--(1+28,12.5);
			\draw [color=red,line width=2pt] (4.5+28,17)--(3+28,17)--(3+28,15)--(1+28,15)--(1+28,25)--(3+28,25)--(3+28,19)--(5+28,19)--(5+28,21)--(6.5+28,21);
			\draw [color=red,line width=2pt] (6.5+28,23)--(5+28,23)--(5+28,25)--(12.5+28,25);
			\draw [color=red,line width=2pt] (17+28,21.5)--(17+28,23)--(15+28,23)--(15+28,25)--(25+28,25)--(25+28,23)--(19+28,23)--(19+28,21)--(21+28,21)--(21+28,19.5);
			\draw [color=red,line width=2pt] (23+28,19.5)--(23+28,21)--(25+28,21)--(25+28,13.5);
			\draw [color=red,line width=2pt] (21.5+28,9)--(23+28,9)--(23+28,11)--(25+28,11)--(25+28,1)--(23+28,1)--(23+28,7)--(21+28,7)--(21+28,5)--(19.5+28,5);
			\draw [color=red,line width=2pt] (19.5+28,3)--(21+28,3)--(21+28,1)--(13.5+28,1);
			\draw [color=violet,line width=1.5pt] (1+28,12.5)--(1+28,13)--(7+28,13)--(7+28,15)--(5+28,15)--(5+28,17)--(4.5+28,17);
			\draw [color=violet,line width=1.5pt] (40.5,25)--(41,25)--(41,19)--(43,19)--(43,21)--(45,21)--(45,21.5);
			\draw [color=brown,line width=1pt] (41,-1)--(41,1)--(41.5,1);
			\draw [color=brown,line width=1pt] (47.5,3)--(47,3)--(47,5)--(47.5,5);
			\draw [color=brown,line width=1pt] (49.5,9)--(49,9)--(49,11)--(47,11)--(47,13)--(51,13)--(51,17)--(49,17)--(49,15)--(47,15)--(47,19)--(45,19)--(45,15)--(43,15)--(43,17)--(41,17)--(41,15)--(39,15)--(39,23)--(37,23)--(37,19)--(35,19)--(35,17)--(37,17)--(37,13)--(45,13)--(45,11)--(39,11)--(39,9)--(37,9)--(37,11)--(31,11)--(31,9)--(35,9)--(35,7)--(41,7)--(41,9)--(47,9)--(47,7)--(43,7)--(43,5)--(45,5)--(45,3)--(41,3)--(41,5)--(37,5)--(37,4.5);
			\draw [color=brown,line width=1pt] (33,6.5)--(33,7)--(31,7)--(31,6.5);
			\draw [color=brown,line width=1pt] (34.5,21)--(35,21)--(35,23)--(34.5,23);
			\draw [color=brown,line width=1pt] (49,19.5)--(49,19)--(51,19)--(51,19.5);
			\draw [color=brown,line width=1pt] (53,13.5)--(53,13)--(55,13);
		\end{tikzpicture}
		\captionof{figure}{Oriental House possible metacell solutions}
	\end{figure}
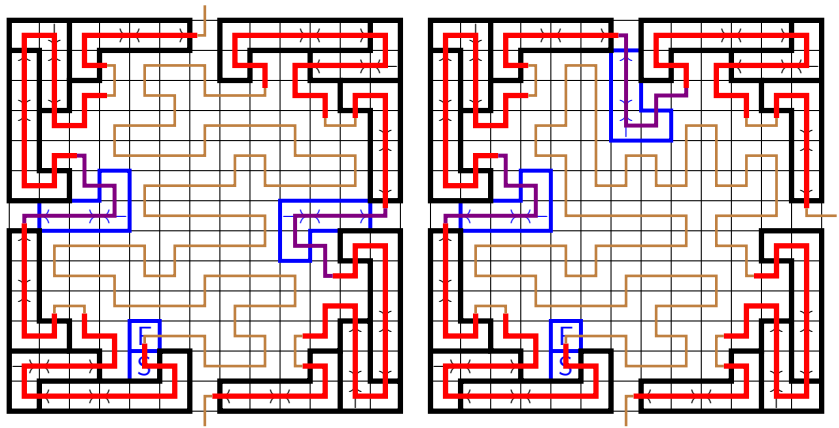
	
	\subsection{Detour: Barred Simple Loop}
	
	Similarly to Oriental House, we notice that the region group depicted in figure 19 has a unique solution, since if the straight lines in the 0-region were vertical then there would necessarily be multiple loops due to the 4-regions. We can extend parts of the 0-region to get a larger forced region, and as depicted in figures 20 and 21 this allows us to make $25\times25$ metacells in a similar manner to Oriental House. Hence, Detour satisfiability is NP-complete as well.
	
	\begin{figure}
		\centering
		\begin{tikzpicture}[scale=0.2,font=\sffamily]
			\foreach \x in {2,4,6,12,14,16} \draw (\x,0)--(\x,12);
			\foreach \y in {2,4,6,8,10} \draw (0,\y)--(8,\y);
			\foreach \y in {2,4,6,8,10} \draw (10,\y)--(18,\y);
			\draw [line width=2pt] (2,0)--(0,0)--(0,12)--(8,12)--(8,0)--(2,0)--(2,4)--(6,4)--(6,0);
			\draw [line width=2pt] (6,12)--(6,8)--(2,8)--(2,12);
			\draw [line width=2pt] (12,0)--(10,0)--(10,12)--(18,12)--(18,0)--(12,0)--(12,4)--(16,4)--(16,0);
			\draw [line width=2pt] (16,12)--(16,8)--(12,8)--(12,12);
			\draw (0.6,11.2) node {\tiny 0};
			\draw (2.6,3.2) node {\tiny 4};
			\draw (2.6,11.2) node {\tiny 4};
			\draw (10.6,11.2) node {\tiny 0};
			\draw (12.6,3.2) node {\tiny 4};
			\draw (12.6,11.2) node {\tiny 4};
			\draw [color=red,line width=2pt] (9,1)--(13,1)--(13,3)--(9,3);
			\draw [color=red,line width=2pt] (19,1)--(15,1)--(15,3)--(19,3);
			\draw [color=red,line width=2pt] (9,5)--(19,5);
			\draw [color=red,line width=2pt] (9,7)--(19,7);
			\draw [color=red,line width=2pt] (9,9)--(13,9)--(13,11)--(9,11);
			\draw [color=red,line width=2pt] (19,9)--(15,9)--(15,11)--(19,11);
		\end{tikzpicture}
		\captionof{figure}{Detour region group with unique solution}
	\end{figure}
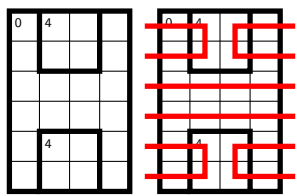
	
	\begin{figure}
		\centering
		\begin{tikzpicture}[scale=0.1,font=\sffamily]
			\foreach \x in {2,4,...,48} \draw (\x,0)--(\x,50);
			\foreach \y in {2,4,...,48} \draw (0,\y)--(50,\y);
			\draw (0,0)--(0,50)--(50,50)--(50,0)--cycle;
			\draw [line width=2pt,color=blue] (24,2)--(26,2)--(26,4)--(24,4); 
			\draw [line width=2pt,color=blue] (2,26)--(2,24)--(4,24)--(4,26);
			\draw [line width=2pt,color=blue] (26,48)--(24,48)--(24,46)--(26,46);
			\draw [line width=2pt,color=blue] (48,24)--(48,26)--(46,26)--(46,24);   
			\draw [line width=2pt] (0,2)--(0,24)--(2,24)--(2,20)--(12,20)--(12,2)--cycle;
			\draw [line width=2pt] (2,0)--(2,12)--(20,12)--(20,2)--(24,2)--(24,4)--(22,4)--(22,0)--cycle;
			\draw [line width=2pt] (14,0)--(14,4)--(18,4)--(18,0);
			\draw [line width=2pt] (14,12)--(14,8)--(18,8)--(18,12);
			\draw [line width=2pt] (20,4)--(22,4);
			\draw [line width=2pt] (0,14)--(4,14)--(4,18)--(0,18);
			\draw [line width=2pt] (12,14)--(8,14)--(8,18)--(12,18);
			\draw [line width=2pt] (0,28)--(0,48)--(12,48)--(12,30)--(2,30)--(2,26)--(4,26)--(4,28)--cycle;
			\draw [line width=2pt] (2,38)--(2,50)--(24,50)--(24,48)--(20,48)--(20,38)--cycle;
			\draw [line width=2pt] (4,28)--(4,30);
			\draw [line width=2pt] (0,32)--(4,32)--(4,36)--(0,36);
			\draw [line width=2pt] (12,32)--(8,32)--(8,36)--(12,36);
			\draw [line width=2pt] (14,50)--(14,46)--(18,46)--(18,50);
			\draw [line width=2pt] (14,38)--(14,42)--(18,42)--(18,38);
			\draw [line width=2pt] (30,38)--(30,48)--(26,48)--(26,46)--(28,46)--(28,50)--(48,50)--(48,38)--cycle;
			\draw [line width=2pt] (38,48)--(50,48)--(50,26)--(48,26)--(48,30)--(38,30)--cycle;
			\draw [line width=2pt] (28,46)--(30,46);
			\draw [line width=2pt] (32,50)--(32,46)--(36,46)--(36,50);
			\draw [line width=2pt] (32,38)--(32,42)--(36,42)--(36,38);
			\draw [line width=2pt] (38,32)--(42,32)--(42,36)--(38,36);
			\draw [line width=2pt] (50,32)--(46,32)--(46,36)--(50,36);
			\draw [line width=2pt] (50,2)--(38,2)--(38,20)--(48,20)--(48,24)--(46,24)--(46,22)--(50,22)--cycle;
			\draw [line width=2pt] (26,0)--(48,0)--(48,12)--(30,12)--(30,2)--(26,2)--cycle;
			\draw [line width=2pt] (46,20)--(46,22);
			\draw [line width=2pt] (32,0)--(32,4)--(36,4)--(36,0);
			\draw [line width=2pt] (32,12)--(32,8)--(36,8)--(36,12);
			\draw [line width=2pt] (38,14)--(42,14)--(42,18)--(38,18);
			\draw [line width=2pt] (50,14)--(46,14)--(46,18)--(50,18);
			\draw (0.8,23) node {\scalebox{0.3}{0}};
			\draw (0.8,17) node {\scalebox{0.3}{4}};
			\draw (8.8,17) node {\scalebox{0.3}{4}};
			\draw (14.8,11) node {\scalebox{0.3}{4}};
			\draw (14.8,3) node {\scalebox{0.3}{4}};
			\draw (12.8,11) node {\scalebox{0.3}{0}};
			\draw (20.8,3) node {\scalebox{0.3}{1}};
			\draw (22.8,3) node {\scalebox{0.3}{0}};
			\draw [color=blue] (24.8,3) node {\scalebox{0.3}{0}};
			\draw (0.8,47) node {\scalebox{0.3}{0}};
			\draw (0.8,35) node {\scalebox{0.3}{4}};
			\draw (8.8,35) node {\scalebox{0.3}{4}};
			\draw (2.8,49) node {\scalebox{0.3}{0}};
			\draw (2.8,29) node {\scalebox{0.3}{1}};
			\draw (2.8,27) node {\scalebox{0.3}{0}};
			\draw [color=blue] (2.8,25) node {\scalebox{0.3}{0}};
			\draw (14.8,49) node {\scalebox{0.3}{4}};
			\draw (14.8,41) node {\scalebox{0.3}{4}};
			\draw (28.8,49) node {\scalebox{0.3}{0}};
			\draw (32.8,49) node {\scalebox{0.3}{4}};
			\draw (32.8,41) node {\scalebox{0.3}{4}};
			\draw (48.8,47) node {\scalebox{0.3}{0}};
			\draw (46.8,35) node {\scalebox{0.3}{4}};
			\draw (38.8,35) node {\scalebox{0.3}{4}};
			\draw [color=blue] (24.8,47) node {\scalebox{0.3}{0}};
			\draw (26.8,47) node {\scalebox{0.3}{0}};
			\draw (28.8,47) node {\scalebox{0.3}{1}};
			\draw [color=blue] (46.8,25) node {\scalebox{0.3}{0}};
			\draw (46.8,23) node {\scalebox{0.3}{0}};
			\draw (46.8,21) node {\scalebox{0.3}{1}};
			\draw (48.8,21) node {\scalebox{0.3}{0}};
			\draw (46.8,17) node {\scalebox{0.3}{4}};
			\draw (38.8,17) node {\scalebox{0.3}{4}};
			\draw (30.8,11) node {\scalebox{0.3}{0}};
			\draw (32.8,11) node {\scalebox{0.3}{4}};
			\draw (32.8,3) node {\scalebox{0.3}{4}};
		\end{tikzpicture}
		\captionof{figure}{Detour metacell}
	\end{figure}
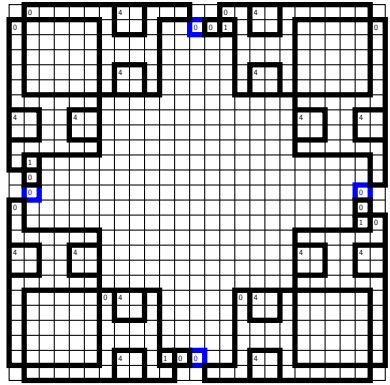

	\begin{figure}
		\centering
		\begin{tikzpicture}[scale=0.1,font=\sffamily]
			\foreach \x in {2,4,...,48} \draw (\x,0)--(\x,50);
			\foreach \y in {2,4,...,48} \draw (0,\y)--(50,\y);
			\draw (0,0)--(0,50)--(50,50)--(50,0)--cycle;
			\draw [line width=2pt,color=blue] (2,26)--(2,24)--(4,24)--(4,26);
			\draw [line width=2pt,color=blue] (48,24)--(48,26)--(46,26)--(46,24);   
			\draw [line width=2pt] (0,2)--(0,24)--(2,24)--(2,20)--(12,20)--(12,2)--cycle;
			\draw [line width=2pt] (2,0)--(2,12)--(20,12)--(20,2)--(24,2)--(24,4)--(22,4)--(22,0)--cycle;
			\draw [line width=2pt] (14,0)--(14,4)--(18,4)--(18,0);
			\draw [line width=2pt] (14,12)--(14,8)--(18,8)--(18,12);
			\draw [line width=2pt] (20,4)--(22,4);
			\draw [line width=2pt] (0,14)--(4,14)--(4,18)--(0,18);
			\draw [line width=2pt] (12,14)--(8,14)--(8,18)--(12,18);
			\draw [line width=2pt] (0,28)--(0,48)--(12,48)--(12,30)--(2,30)--(2,26)--(4,26)--(4,28)--cycle;
			\draw [line width=2pt] (2,38)--(2,50)--(24,50)--(24,48)--(20,48)--(20,38)--cycle;
			\draw [line width=2pt] (4,28)--(4,30);
			\draw [line width=2pt] (0,32)--(4,32)--(4,36)--(0,36);
			\draw [line width=2pt] (12,32)--(8,32)--(8,36)--(12,36);
			\draw [line width=2pt] (14,50)--(14,46)--(18,46)--(18,50);
			\draw [line width=2pt] (14,38)--(14,42)--(18,42)--(18,38);
			\draw [line width=2pt] (30,38)--(30,48)--(26,48)--(26,46)--(28,46)--(28,50)--(48,50)--(48,38)--cycle;
			\draw [line width=2pt] (38,48)--(50,48)--(50,26)--(48,26)--(48,30)--(38,30)--cycle;
			\draw [line width=2pt] (28,46)--(30,46);
			\draw [line width=2pt] (32,50)--(32,46)--(36,46)--(36,50);
			\draw [line width=2pt] (32,38)--(32,42)--(36,42)--(36,38);
			\draw [line width=2pt] (38,32)--(42,32)--(42,36)--(38,36);
			\draw [line width=2pt] (50,32)--(46,32)--(46,36)--(50,36);
			\draw [line width=2pt] (50,2)--(38,2)--(38,20)--(48,20)--(48,24)--(46,24)--(46,22)--(50,22)--cycle;
			\draw [line width=2pt] (26,0)--(48,0)--(48,12)--(30,12)--(30,2)--(26,2)--cycle;
			\draw [line width=2pt] (46,20)--(46,22);
			\draw [line width=2pt] (32,0)--(32,4)--(36,4)--(36,0);
			\draw [line width=2pt] (32,12)--(32,8)--(36,8)--(36,12);
			\draw [line width=2pt] (38,14)--(42,14)--(42,18)--(38,18);
			\draw [line width=2pt] (50,14)--(46,14)--(46,18)--(50,18);
			\draw (0.8,23) node {\scalebox{0.3}{0}};
			\draw (0.8,17) node {\scalebox{0.3}{4}};
			\draw (8.8,17) node {\scalebox{0.3}{4}};
			\draw (14.8,11) node {\scalebox{0.3}{4}};
			\draw (14.8,3) node {\scalebox{0.3}{4}};
			\draw (12.8,11) node {\scalebox{0.3}{0}};
			\draw (20.8,3) node {\scalebox{0.3}{1}};
			\draw (22.8,3) node {\scalebox{0.3}{0}};
			\draw (0.8,47) node {\scalebox{0.3}{0}};
			\draw (0.8,35) node {\scalebox{0.3}{4}};
			\draw (8.8,35) node {\scalebox{0.3}{4}};
			\draw (2.8,49) node {\scalebox{0.3}{0}};
			\draw (2.8,29) node {\scalebox{0.3}{1}};
			\draw (2.8,27) node {\scalebox{0.3}{0}};
			\draw [color=blue] (2.8,25) node {\scalebox{0.3}{0}};
			\draw (14.8,49) node {\scalebox{0.3}{4}};
			\draw (14.8,41) node {\scalebox{0.3}{4}};
			\draw (28.8,49) node {\scalebox{0.3}{0}};
			\draw (32.8,49) node {\scalebox{0.3}{4}};
			\draw (32.8,41) node {\scalebox{0.3}{4}};
			\draw (48.8,47) node {\scalebox{0.3}{0}};
			\draw (46.8,35) node {\scalebox{0.3}{4}};
			\draw (38.8,35) node {\scalebox{0.3}{4}};
			\draw (26.8,47) node {\scalebox{0.3}{0}};
			\draw (28.8,47) node {\scalebox{0.3}{1}};
			\draw [color=blue] (46.8,25) node {\scalebox{0.3}{0}};
			\draw (46.8,23) node {\scalebox{0.3}{0}};
			\draw (46.8,21) node {\scalebox{0.3}{1}};
			\draw (48.8,21) node {\scalebox{0.3}{0}};
			\draw (46.8,17) node {\scalebox{0.3}{4}};
			\draw (38.8,17) node {\scalebox{0.3}{4}};
			\draw (30.8,11) node {\scalebox{0.3}{0}};
			\draw (32.8,11) node {\scalebox{0.3}{4}};
			\draw (32.8,3) node {\scalebox{0.3}{4}};
			
			\draw [color=brown,line width=1pt] (25,-1)--(25,1)--(25.5,1);
			\draw [color=red,line width=2pt] (25.5,1)--(33,1)--(33,3)--(29.5,3);
			\draw [color=brown,line width=1pt] (29.5,3)--(25,3)--(25,5)--(29.5,5);
			\draw [color=red,line width=2pt] (29.5,5)--(38.5,5);
			\draw [color=brown,line width=1pt] (38.5,5)--(39,5)--(39,3)--(38.5,3);
			\draw [color=red,line width=2pt] (38.5,3)--(35,3)--(35,1)--(49,1)--(49,15)--(47,15)--(47,11.5);
			\draw [color=brown,line width=1pt] (47,11.5)--(47,4.5);
			\draw [color=red,line width=2pt] (47,4.5)--(47,3)--(45.5,3);
			\draw [color=brown,line width=1pt] (45.5,3)--(45,3)--(45,11.5);
			\draw [color=red,line width=2pt] (45,11.5)--(45,21)--(47,21)--(47,17)--(49,17)--(49,23)--(45.5,23);
			\draw [color=brown,line width=1pt] (45.5,23)--(43,23)--(43,20.5);
			\draw [color=red,line width=2pt] (43,20.5)--(43,11.5);
			\draw [color=brown,line width=1pt] (43,11.5)--(43,3)--(41,3)--(41,11.5);
			\draw [color=red,line width=2pt] (41,11.5)--(41,15)--(39,15)--(39,11)--(35,11)--(35,9)--(38.5,9);
			\draw [color=brown,line width=1pt] (38.5,9)--(39,9)--(39,7)--(38.5,7);
			\draw [color=red,line width=2pt] (38.5,7)--(29.5,7);
			\draw [color=brown,line width=1pt] (29.5,7)--(25,7)--(25,9)--(29.5,9);
			\draw [color=red,line width=2pt] (29.5,9)--(33,9)--(33,11)--(29.5,11);
			\draw [color=brown,line width=1pt] (29.5,11)--(23,11)--(23,9)--(20.5,9);
			\draw [color=red,line width=2pt] (20.5,9)--(17,9)--(17,11)--(20.5,11);
			\draw [color=brown,line width=1pt] (20.5,11)--(21,11)--(21,13)--(35.5,13);
			\draw [color=red,line width=2pt] (35.5,13)--(37,13)--(37,14.5);
			\draw [color=brown,line width=1pt] (37,14.5)--(37,15)--(19,15)--(19,13)--(17,13)--(17,15)--(15,15)--(15,13)--(14.5,13);
			\draw [color=red,line width=2pt] (14.5,13)--(13,13)--(13,14.5);
			\draw [color=brown,line width=1pt] (13,14.5)--(13,17)--(37,17)--(37,19)--(13,19)--(13,21)--(39,21)--(39,20.5);
			\draw [color=red,line width=2pt] (39,20.5)--(39,17)--(41,17)--(41,20.5);
			\draw [color=brown,line width=1pt] (41,20.5)--(41,23)--(11,23)--(11,20.5);
			\draw [color=red,line width=2pt] (11,20.5)--(11,17)--(9,17)--(9,20.5);
			\draw [color=brown,line width=1pt] (9,20.5)--(9,23)--(7,23)--(7,20.5);
			\draw [color=red,line width=2pt] (7,20.5)--(7,11.5);
			\draw [color=brown,line width=1pt] (7,11.5)--(7,11)--(9,11)--(9,11.5);
			\draw [color=red,line width=2pt] (9,11.5)--(9,15)--(11,15)--(11,11)--(15,11)--(15,9)--(11.5,9);
			\draw [color=brown,line width=1pt] (11.5,9)--(3,9)--(3,7)--(11.5,7);
			\draw [color=red,line width=2pt] (11.5,7)--(20.5,7);
			\draw [color=brown,line width=1pt] (20.5,7)--(23,7)--(23,4.5);
			\draw [color=red,line width=2pt] (23,4.5)--(23,1)--(17,1)--(17,3)--(21,3)--(21,5)--(11.5,5);
			\draw [color=brown,line width=1pt] (11.5,5)--(3,5)--(3,4.5);
			\draw [color=red,line width=2pt] (3,4.5)--(3,3)--(4.5,3);
			\draw [color=brown,line width=1pt] (4.5,3)--(11.5,3);
			\draw [color=red,line width=2pt] (11.5,3)--(15,3)--(15,1)--(1,1)--(1,15)--(3,15)--(3,11.5);
			\draw [color=brown,line width=1pt] (3,11.5)--(3,11)--(5,11)--(5,11.5);
			\draw [color=red,line width=2pt] (5,11.5)--(5,20.5);
			\draw [color=brown,line width=1pt] (5,20.5)--(5,23)--(4.5,23);
			\draw [color=violet,line width=1.5pt] (4.5,23)--(3,23)--(3,20.5);
			\draw [color=red,line width=2pt] (3,20.5)--(3,17)--(1,17)--(1,24.5);
			\draw [color=violet,line width=1.5pt] (1,24.5)--(1,25)--(4.5,25);
			\draw [color=red,line width=2pt] (49,25.5)--(49,33)--(47,33)--(47,29.5);
			\draw [color=red,line width=2pt] (45,29.5)--(45,38.5);
			\draw [color=brown,line width=1pt] (45,38.5)--(45,39)--(47,39)--(47,38.5);
			\draw [color=red,line width=2pt] (47,38.5)--(47,35)--(49,35)--(49,49)--(35,49)--(35,47)--(38.5,47);
			\draw [color=brown,line width=1pt] (38.5,47)--(45.5,47);
			\draw [color=red,line width=2pt] (45.5,47)--(47,47)--(47,45.5);
			\draw [color=brown,line width=1pt] (47,45.5)--(47,45)--(38.5,45);
			\draw [color=red,line width=2pt] (38.5,45)--(29,45)--(29,47)--(33,47)--(33,49)--(27,49)--(27,45.5);
			\draw [color=brown,line width=1pt] (27,45.5)--(27,43)--(29.5,43);
			\draw [color=red,line width=2pt] (29.5,43)--(38.5,43);
			\draw [color=brown,line width=1pt] (38.5,43)--(47,43)--(47,41)--(38.5,41);
			\draw [color=red,line width=2pt] (38.5,41)--(35,41)--(35,39)--(39,39)--(39,35)--(41,35)--(41,38.5);
			\draw [color=brown,line width=1pt] (41,38.5)--(41,39)--(43,39)--(43,38.5);
			\draw [color=red,line width=2pt] (43,38.5)--(43,29.5);
			\draw [color=red,line width=2pt] (41,29.5)--(41,33)--(39,33)--(39,29.5);
			\draw [color=red,line width=2pt] (9,29.5)--(9,33)--(11,33)--(11,29.5);
			\draw [color=red,line width=2pt] (37,35.5)--(37,37)--(35.5,37);
			\draw [color=red,line width=2pt] (13,35.5)--(13,37)--(14.5,37);
			\draw [color=red,line width=2pt] (29.5,39)--(33,39)--(33,41)--(29.5,41);
			\draw [color=red,line width=2pt] (20.5,39)--(17,39)--(17,41)--(20.5,41);
			\draw [color=red,line width=2pt] (20.5,43)--(11.5,43);
			\draw [color=brown,line width=1pt] (11.5,43)--(11,43)--(11,41)--(11.5,41);
			\draw [color=red,line width=2pt] (11.5,41)--(15,41)--(15,39)--(11,39)--(11,35)--(9,35)--(9,38.5);
			\draw [color=brown,line width=1pt] (9,38.5)--(9,47)--(7,47)--(7,38.5);
			\draw [color=red,line width=2pt] (7,38.5)--(7,29.5);
			\draw [color=brown,line width=1pt] (7,29.5)--(7,27)--(4.5,27);
			\draw [color=red,line width=2pt] (4.5,27)--(1,27)--(1,33)--(3,33)--(3,29)--(5,29)--(5,38.5);
			\draw [color=brown,line width=1pt] (5,38.5)--(5,47)--(4.5,47);
			\draw [color=red,line width=2pt] (4.5,47)--(3,47)--(3,45.5);
			\draw [color=brown,line width=1pt] (3,45.5)--(3,38.5);
			\draw [color=red,line width=2pt] (3,38.5)--(3,35)--(1,35)--(1,49)--(15,49)--(15,47)--(11.5,47);
			\draw [color=brown,line width=1pt] (11.5,47)--(11,47)--(11,45)--(11.5,45);
			\draw [color=red,line width=2pt] (11.5,45)--(20.5,45);
			\draw [color=red,line width=2pt] (20.5,47)--(17,47)--(17,49)--(24.5,49);
			
			\draw [color=brown,line width=1pt] (4.5,25)--(45.5,25);
			\draw [color=violet,line width=1.5pt] (45.5,25)--(49,25)--(49,25.5);
			\draw [color=violet,line width=1.5pt] (47,29.5)--(47,27)--(45.5,27);
			\draw [color=brown,line width=1pt] (45.5,27)--(45,27)--(45,29.5);
			\draw [color=brown,line width=1pt] (43,29.5)--(43,27)--(41,27)--(41,29.5);
			\draw [color=brown,line width=1pt] (39,29.5)--(39,27)--(9,27)--(9,29.5);
			\draw [color=brown,line width=1pt] (11,29.5)--(11,29)--(37,29)--(37,31)--(13,31)--(13,33)--(37,33)--(37,35.5);
			\draw [color=brown,line width=1pt] (35.5,37)--(35,37)--(35,35)--(33,35)--(33,37)--(31,37)--(31,35)--(13,35)--(13,35.5);
			\draw [color=brown,line width=1pt] (14.5,37)--(29,37)--(29,39)--(29.5,39);
			\draw [color=brown,line width=1pt] (29.5,41)--(27,41)--(27,39)--(20.5,39);
			\draw [color=brown,line width=1pt] (20.5,41)--(25,41)--(25,43)--(20.5,43);
			\draw [color=brown,line width=1pt] (20.5,45)--(25,45)--(25,47)--(20.5,47);
			\draw [color=brown,line width=1pt] (24.5,49)--(25,49)--(25,51);		
	\end{tikzpicture}
	\end{figure}

	\begin{figure}
		\centering
		\begin{tikzpicture}[scale=0.1,font=\sffamily]
			\foreach \x in {2,4,...,48} \draw (\x,0)--(\x,50);
			\foreach \y in {2,4,...,48} \draw (0,\y)--(50,\y);
			\draw (0,0)--(0,50)--(50,50)--(50,0)--cycle;
			\draw [line width=2pt,color=blue] (2,26)--(2,24)--(4,24)--(4,26);
			\draw [line width=2pt,color=blue] (26,48)--(24,48)--(24,46)--(26,46);  
			\draw [line width=2pt] (0,2)--(0,24)--(2,24)--(2,20)--(12,20)--(12,2)--cycle;
			\draw [line width=2pt] (2,0)--(2,12)--(20,12)--(20,2)--(24,2)--(24,4)--(22,4)--(22,0)--cycle;
			\draw [line width=2pt] (14,0)--(14,4)--(18,4)--(18,0);
			\draw [line width=2pt] (14,12)--(14,8)--(18,8)--(18,12);
			\draw [line width=2pt] (20,4)--(22,4);
			\draw [line width=2pt] (0,14)--(4,14)--(4,18)--(0,18);
			\draw [line width=2pt] (12,14)--(8,14)--(8,18)--(12,18);
			\draw [line width=2pt] (0,28)--(0,48)--(12,48)--(12,30)--(2,30)--(2,26)--(4,26)--(4,28)--cycle;
			\draw [line width=2pt] (2,38)--(2,50)--(24,50)--(24,48)--(20,48)--(20,38)--cycle;
			\draw [line width=2pt] (4,28)--(4,30);
			\draw [line width=2pt] (0,32)--(4,32)--(4,36)--(0,36);
			\draw [line width=2pt] (12,32)--(8,32)--(8,36)--(12,36);
			\draw [line width=2pt] (14,50)--(14,46)--(18,46)--(18,50);
			\draw [line width=2pt] (14,38)--(14,42)--(18,42)--(18,38);
			\draw [line width=2pt] (30,38)--(30,48)--(26,48)--(26,46)--(28,46)--(28,50)--(48,50)--(48,38)--cycle;
			\draw [line width=2pt] (38,48)--(50,48)--(50,26)--(48,26)--(48,30)--(38,30)--cycle;
			\draw [line width=2pt] (28,46)--(30,46);
			\draw [line width=2pt] (32,50)--(32,46)--(36,46)--(36,50);
			\draw [line width=2pt] (32,38)--(32,42)--(36,42)--(36,38);
			\draw [line width=2pt] (38,32)--(42,32)--(42,36)--(38,36);
			\draw [line width=2pt] (50,32)--(46,32)--(46,36)--(50,36);
			\draw [line width=2pt] (50,2)--(38,2)--(38,20)--(48,20)--(48,24)--(46,24)--(46,22)--(50,22)--cycle;
			\draw [line width=2pt] (26,0)--(48,0)--(48,12)--(30,12)--(30,2)--(26,2)--cycle;
			\draw [line width=2pt] (46,20)--(46,22);
			\draw [line width=2pt] (32,0)--(32,4)--(36,4)--(36,0);
			\draw [line width=2pt] (32,12)--(32,8)--(36,8)--(36,12);
			\draw [line width=2pt] (38,14)--(42,14)--(42,18)--(38,18);
			\draw [line width=2pt] (50,14)--(46,14)--(46,18)--(50,18);
			\draw (0.8,23) node {\scalebox{0.3}{0}};
			\draw (0.8,17) node {\scalebox{0.3}{4}};
			\draw (8.8,17) node {\scalebox{0.3}{4}};
			\draw (14.8,11) node {\scalebox{0.3}{4}};
			\draw (14.8,3) node {\scalebox{0.3}{4}};
			\draw (12.8,11) node {\scalebox{0.3}{0}};
			\draw (20.8,3) node {\scalebox{0.3}{1}};
			\draw (22.8,3) node {\scalebox{0.3}{0}};
			\draw (0.8,47) node {\scalebox{0.3}{0}};
			\draw (0.8,35) node {\scalebox{0.3}{4}};
			\draw (8.8,35) node {\scalebox{0.3}{4}};
			\draw (2.8,49) node {\scalebox{0.3}{0}};
			\draw (2.8,29) node {\scalebox{0.3}{1}};
			\draw (2.8,27) node {\scalebox{0.3}{0}};
			\draw [color=blue] (2.8,25) node {\scalebox{0.3}{0}};
			\draw (14.8,49) node {\scalebox{0.3}{4}};
			\draw (14.8,41) node {\scalebox{0.3}{4}};
			\draw (28.8,49) node {\scalebox{0.3}{0}};
			\draw (32.8,49) node {\scalebox{0.3}{4}};
			\draw (32.8,41) node {\scalebox{0.3}{4}};
			\draw (48.8,47) node {\scalebox{0.3}{0}};
			\draw (46.8,35) node {\scalebox{0.3}{4}};
			\draw (38.8,35) node {\scalebox{0.3}{4}};
			\draw [color=blue] (24.8,47) node {\scalebox{0.3}{0}};
			\draw (26.8,47) node {\scalebox{0.3}{0}};
			\draw (28.8,47) node {\scalebox{0.3}{1}};
			\draw (46.8,23) node {\scalebox{0.3}{0}};
			\draw (46.8,21) node {\scalebox{0.3}{1}};
			\draw (48.8,21) node {\scalebox{0.3}{0}};
			\draw (46.8,17) node {\scalebox{0.3}{4}};
			\draw (38.8,17) node {\scalebox{0.3}{4}};
			\draw (30.8,11) node {\scalebox{0.3}{0}};
			\draw (32.8,11) node {\scalebox{0.3}{4}};
			\draw (32.8,3) node {\scalebox{0.3}{4}};
			
			\draw [color=brown,line width=1pt] (25,-1)--(25,1)--(25.5,1);
			\draw [color=red,line width=2pt] (25.5,1)--(33,1)--(33,3)--(29.5,3);
			\draw [color=brown,line width=1pt] (29.5,3)--(25,3)--(25,5)--(29.5,5);
			\draw [color=red,line width=2pt] (29.5,5)--(38.5,5);
			\draw [color=brown,line width=1pt] (38.5,5)--(39,5)--(39,3)--(38.5,3);
			\draw [color=red,line width=2pt] (38.5,3)--(35,3)--(35,1)--(49,1)--(49,15)--(47,15)--(47,11.5);
			\draw [color=brown,line width=1pt] (47,11.5)--(47,4.5);
			\draw [color=red,line width=2pt] (47,4.5)--(47,3)--(45.5,3);
			\draw [color=brown,line width=1pt] (45.5,3)--(45,3)--(45,11.5);
			\draw [color=red,line width=2pt] (45,11.5)--(45,21)--(47,21)--(47,17)--(49,17)--(49,23)--(45.5,23);
			\draw [color=brown,line width=1pt] (45.5,23)--(43,23)--(43,20.5);
			\draw [color=red,line width=2pt] (43,20.5)--(43,11.5);
			\draw [color=brown,line width=1pt] (43,11.5)--(43,3)--(41,3)--(41,11.5);
			\draw [color=red,line width=2pt] (41,11.5)--(41,15)--(39,15)--(39,11)--(35,11)--(35,9)--(38.5,9);
			\draw [color=brown,line width=1pt] (38.5,9)--(39,9)--(39,7)--(38.5,7);
			\draw [color=red,line width=2pt] (38.5,7)--(29.5,7);
			\draw [color=brown,line width=1pt] (29.5,7)--(25,7)--(25,9)--(29.5,9);
			\draw [color=red,line width=2pt] (29.5,9)--(33,9)--(33,11)--(29.5,11);
			\draw [color=brown,line width=1pt] (29.5,11)--(23,11)--(23,9)--(20.5,9);
			\draw [color=red,line width=2pt] (20.5,9)--(17,9)--(17,11)--(20.5,11);
			\draw [color=brown,line width=1pt] (20.5,11)--(21,11)--(21,13)--(35.5,13);
			\draw [color=red,line width=2pt] (35.5,13)--(37,13)--(37,14.5);
			\draw [color=brown,line width=1pt] (37,14.5)--(37,15)--(19,15)--(19,13)--(17,13)--(17,15)--(15,15)--(15,13)--(14.5,13);
			\draw [color=red,line width=2pt] (14.5,13)--(13,13)--(13,14.5);
			\draw [color=brown,line width=1pt] (13,14.5)--(13,17)--(37,17)--(37,19)--(13,19)--(13,21)--(39,21)--(39,20.5);
			\draw [color=red,line width=2pt] (39,20.5)--(39,17)--(41,17)--(41,20.5);
			\draw [color=brown,line width=1pt] (41,20.5)--(41,23)--(11,23)--(11,20.5);
			\draw [color=red,line width=2pt] (11,20.5)--(11,17)--(9,17)--(9,20.5);
			\draw [color=brown,line width=1pt] (9,20.5)--(9,23)--(7,23)--(7,20.5);
			\draw [color=red,line width=2pt] (7,20.5)--(7,11.5);
			\draw [color=brown,line width=1pt] (7,11.5)--(7,11)--(9,11)--(9,11.5);
			\draw [color=red,line width=2pt] (9,11.5)--(9,15)--(11,15)--(11,11)--(15,11)--(15,9)--(11.5,9);
			\draw [color=brown,line width=1pt] (11.5,9)--(3,9)--(3,7)--(11.5,7);
			\draw [color=red,line width=2pt] (11.5,7)--(20.5,7);
			\draw [color=brown,line width=1pt] (20.5,7)--(23,7)--(23,4.5);
			\draw [color=red,line width=2pt] (23,4.5)--(23,1)--(17,1)--(17,3)--(21,3)--(21,5)--(11.5,5);
			\draw [color=brown,line width=1pt] (11.5,5)--(3,5)--(3,4.5);
			\draw [color=red,line width=2pt] (3,4.5)--(3,3)--(4.5,3);
			\draw [color=brown,line width=1pt] (4.5,3)--(11.5,3);
			\draw [color=red,line width=2pt] (11.5,3)--(15,3)--(15,1)--(1,1)--(1,15)--(3,15)--(3,11.5);
			\draw [color=brown,line width=1pt] (3,11.5)--(3,11)--(5,11)--(5,11.5);
			\draw [color=red,line width=2pt] (5,11.5)--(5,20.5);
			\draw [color=brown,line width=1pt] (5,20.5)--(5,23)--(4.5,23);
			\draw [color=violet,line width=1.5pt] (4.5,23)--(3,23)--(3,20.5);
			\draw [color=red,line width=2pt] (3,20.5)--(3,17)--(1,17)--(1,24.5);
			\draw [color=violet,line width=1.5pt] (1,24.5)--(1,25)--(4.5,25);
			\draw [color=red,line width=2pt] (49,25.5)--(49,33)--(47,33)--(47,29.5);
			\draw [color=red,line width=2pt] (45,29.5)--(45,38.5);
			\draw [color=brown,line width=1pt] (45,38.5)--(45,39)--(47,39)--(47,38.5);
			\draw [color=red,line width=2pt] (47,38.5)--(47,35)--(49,35)--(49,49)--(35,49)--(35,47)--(38.5,47);
			\draw [color=brown,line width=1pt] (38.5,47)--(45.5,47);
			\draw [color=red,line width=2pt] (45.5,47)--(47,47)--(47,45.5);
			\draw [color=brown,line width=1pt] (47,45.5)--(47,45)--(38.5,45);
			\draw [color=red,line width=2pt] (38.5,45)--(29,45)--(29,47)--(33,47)--(33,49)--(27,49)--(27,45.5);
			\draw [color=brown,line width=1pt] (27,45.5)--(27,43)--(29.5,43);
			\draw [color=red,line width=2pt] (29.5,43)--(38.5,43);
			\draw [color=brown,line width=1pt] (38.5,43)--(47,43)--(47,41)--(38.5,41);
			\draw [color=red,line width=2pt] (38.5,41)--(35,41)--(35,39)--(39,39)--(39,35)--(41,35)--(41,38.5);
			\draw [color=brown,line width=1pt] (41,38.5)--(41,39)--(43,39)--(43,38.5);
			\draw [color=red,line width=2pt] (43,38.5)--(43,29.5);
			\draw [color=red,line width=2pt] (41,29.5)--(41,33)--(39,33)--(39,29.5);
			\draw [color=red,line width=2pt] (9,29.5)--(9,33)--(11,33)--(11,29.5);
			\draw [color=red,line width=2pt] (37,35.5)--(37,37)--(35.5,37);
			\draw [color=red,line width=2pt] (13,35.5)--(13,37)--(14.5,37);
			\draw [color=red,line width=2pt] (29.5,39)--(33,39)--(33,41)--(29.5,41);
			\draw [color=red,line width=2pt] (20.5,39)--(17,39)--(17,41)--(20.5,41);
			\draw [color=red,line width=2pt] (20.5,43)--(11.5,43);
			\draw [color=brown,line width=1pt] (11.5,43)--(11,43)--(11,41)--(11.5,41);
			\draw [color=red,line width=2pt] (11.5,41)--(15,41)--(15,39)--(11,39)--(11,35)--(9,35)--(9,38.5);
			\draw [color=brown,line width=1pt] (9,38.5)--(9,47)--(7,47)--(7,38.5);
			\draw [color=red,line width=2pt] (7,38.5)--(7,29.5);
			\draw [color=brown,line width=1pt] (7,29.5)--(7,27)--(4.5,27);
			\draw [color=red,line width=2pt] (4.5,27)--(1,27)--(1,33)--(3,33)--(3,29)--(5,29)--(5,38.5);
			\draw [color=brown,line width=1pt] (5,38.5)--(5,47)--(4.5,47);
			\draw [color=red,line width=2pt] (4.5,47)--(3,47)--(3,45.5);
			\draw [color=brown,line width=1pt] (3,45.5)--(3,38.5);
			\draw [color=red,line width=2pt] (3,38.5)--(3,35)--(1,35)--(1,49)--(15,49)--(15,47)--(11.5,47);
			\draw [color=brown,line width=1pt] (11.5,47)--(11,47)--(11,45)--(11.5,45);
			\draw [color=red,line width=2pt] (11.5,45)--(20.5,45);
			\draw [color=red,line width=2pt] (20.5,47)--(17,47)--(17,49)--(24.5,49);
			
			\draw [color=brown,line width=1pt] (4.5,25)--(9,25)--(9,29.5);
			\draw [color=brown,line width=1pt] (11,29.5)--(11,25)--(13,25)--(13,35.5);
			\draw [color=brown,line width=1pt] (14.5,37)--(15,37)--(15,25)--(17,25)--(17,37)--(19,37)--(19,25)--(21,25)--(21,39)--(20.5,39);
			\draw [color=brown,line width=1pt] (20.5,41)--(21,41)--(21,43)--(20.5,43);
			\draw [color=brown,line width=1pt] (20.5,45)--(23,45)--(23,45.5);
			\draw [color=violet,line width=1.5pt] (23,45.5)--(23,47)--(20.5,47);
			\draw [color=violet,line width=1.5pt] (24.5,49)--(25,49)--(25,45.5);
			\draw [color=brown,line width=1pt] (25,45.5)--(25,43)--(23,43)--(23,25)--(25,25)--(25,41)--(29.5,41)--(29.5,39)--(27,39)--(27,25)--(29,25)--(29,37)--(31,37)--(31,25)--(33,25)--(33,37)--(35.5,37);
			\draw [color=brown,line width=1pt] (37,35.5)--(37,35)--(35,35)--(35,33)--(37,33)--(37,31)--(35,31)--(35,29)--(37,29)--(37,27)--(35,27)--(35,25)--(39,25)--(39,29.5);
			\draw [color=brown,line width=1pt] (41,29.5)--(41,25)--(43,25)--(43,29.5);
			\draw [color=brown,line width=1pt] (45,29.5)--(45,25)--(47,25)--(47,29.5);
			\draw [color=brown,line width=1pt] (49,25.5)--(49,25)--(51,25);
		\end{tikzpicture}
		\captionof{figure}{Detour possible metacell solutions}
	\end{figure}

	\section{Concluding Remarks and Acknowledgements}
	
	The techniques presented in this paper are quite likely applicable to a wide variety of other loop puzzles. For example, we can achieve a reduction from the Haisu reduction to Country Road \cite{puzll} by setting the numbers in each country road region to the number of cells in the region: this ensures each region is entered and exited once. Similarly, due to the use of all-turns or all-straights regions, the Detour construction immediately reduces to Dutch Loop \cite{puzll}, although a metacell construction is probably achievable in a wide variety of other genres.
	
	Thanks to William Hu for testsolving the example puzzles, and to both William Hu and Guowen Zhang for the invention of these genres. Also thanks to Charles Li for helping to proofread the paper.

\end{document}